\documentclass[a4paper, 11 pt]{article}
\usepackage{a4wide, amsmath, amssymb, mathtools, yfonts}
\usepackage{comment}
\usepackage{tikz}
\usepackage{tikz-cd}
\usepackage[all]{xy}
\usepackage[utf8]{inputenc}
\usepackage{amsthm, mathrsfs}
\usepackage[english]{babel}
\usepackage{hyperref}
\usepackage{authblk}
\usepackage[OT2,T1]{fontenc}
\DeclareSymbolFont{cyrletters}{OT2}{wncyr}{m}{n}
\DeclareMathSymbol{\Sha}{\mathalpha}{cyrletters}{"58}

\numberwithin{equation}{section}
\newtheorem{lemma}{Lemma}[section]
\newtheorem{theorem}[lemma]{Theorem}
\newtheorem{proposition}[lemma]{Proposition}
\newtheorem{corollary}[lemma]{Corollary}

\newtheorem{remark}[lemma]{Remark}

\theoremstyle{definition}
\newtheorem{mydef}[lemma]{Definition}

\newcommand{\Z}{\mathbb{Z}}
\newcommand{\Q}{\mathbb{Q}}
\newcommand{\C}{\mathbb{C}}

\newcommand{\R}{\mathbb{R}}
\newcommand{\FF}{\mathbb{F}}

\newcommand{\Hom}{\mathrm{Hom}}
\newcommand{\Frob}{\textup{Frob}}

\newcommand\Vol{\mathrm{Vol}}
\newcommand\Gal{\mathrm{Gal}}

\newcommand{\res}{\textup{res}}

\usepackage{epigraph}
\allowdisplaybreaks

\title{\vspace{-\baselineskip}\sffamily\bfseries Hilbert's tenth problem via additive combinatorics}
\author[1]{Peter Koymans\thanks{Mathematisch Instituut, Universiteit Utrecht, Postbus 80.010, 3508 TA Utrecht, The Netherlands, p.h.koymans@uu.nl}}
\author[2]{Carlo Pagano\thanks{Department of Mathematics and Statistics, Montreal, Quebec H3G 1M8, Canada, carlein90@gmail.com}}
\affil[1]{Utrecht University}
\affil[2]{Concordia University}
\date{\today}

\begin{document}
\maketitle

\begin{abstract}
For all infinite rings $R$ that are finitely generated over $\Z$, we show that Hilbert's tenth problem has a negative answer. This is accomplished by constructing elliptic curves $E$ without rank growth in certain quadratic extensions $L/K$. To achieve such a result unconditionally, our key innovation is to use elliptic curves $E$ with full rational $2$-torsion which allows us to combine techniques from additive combinatorics with $2$-descent.
\end{abstract}

\epigraph{\emph{Wir m\"ussen wissen \\
Wir werden wissen}}{David Hilbert}

\section{Introduction}
\subsection{History of the problem}
At the 1900 mathematical congress in Paris, Hilbert introduced his famous list of 23 problems. The tenth problem of his list is nowadays interpreted as asking for an algorithm to decide the solvability of Diophantine equations. G\"odel dealt a major blow to this program in 1931 with his incompleteness theorems, although his results certainly did not rule out the existence of such an algorithm.

In 1950, Julia Robinson proved that exponentiation is Diophantine conditional under a hypothesis that later became known as ``J.R.''. One major ingredient in this step is the basic theory of Pell's equation. Not much later, Davis, Putnam and Robinson \cite{DPR} showed that the hypothesis J.R.~implies that every recursively enumerable set is Diophantine, hence Hilbert's tenth problem is unsolvable over the integers. Finally, Matiyasevich \cite{Ma} proved J.R.~using properties of Fibonacci numbers, thus settling Hilbert's tenth problem in full. Because of the combined efforts of these four authors, this result became known as the MRDP theorem.

Nevertheless, the field remained an active research area as the MRDP theorem only settles the case of integer polynomials and integer solutions, while it is natural to wonder what happens for other number-theoretic rings. One of the main difficulties was to find an appropriate substitute for Pell's equation to show that exponentiation is Diophantine. Fairly early on researchers realized that algebraic groups of rank $1$ could fulfill such a role, but rigorously establishing their existence has been a difficult problem. As it stands, results from the 1980s give a negative answer to Hilbert's tenth problem for the ring of integers of the following number fields $K$:
\begin{itemize}
\item $K$ is totally real or a quadratic extension of a totally real field \cite{DL, Denef},
\item $[K : \Q] = 4$ and $K$ is not totally real \cite{DL},
\item $K$ has exactly one complex place \cite{Ph, Sh, Vi},
\item $K$ is contained in a number field $L$ such that $\Z$ is Diophantine over $O_L$ \cite{SS}.
\end{itemize}
We also mention the recent works of Garcia-Fritz--Pasten \cite{GP}, Kundu--Lei--Sprung \cite{KLS} and Shnidman--Weiss \cite{Shnidman--Weiss}, who give a negative answer for a large number of sextic fields. There are also related results on Diophantine stability in \cite{MRL} and \cite{RW}.

\subsection{Our results}
After the completion of the MRDP theorem, elliptic curves became a particularly natural candidate to replace the role of Pell's equation. This was briefly suggested by Denef in \cite{Denef}. A rigorous connection was first established by Poonen \cite{Poonen} and independently by Cornelissen--Pheidas--Zahidi \cite{CPZ}. These results were further strengthened by Shlapentokh \cite{Shl} and Mazur--Rubin--Shlapentokh \cite{MRS}. 

However, these results came with one major downside: unlike Pell's equation, it is not easy to show that there exists an elliptic curve of rank $1$ over a given number field. In fact, only after the release of our initial preprint, this was shown to be true independently by Zywina \cite{Zywina} and the authors \cite{KPE}.

However, assuming finiteness of the Tate--Shafarevich group, Mazur--Rubin \cite{MR} were able to use the aforementioned results to show that Hilbert's tenth problem is undecidable for any infinite ring $R$ that is finitely generated over $\Z$. Further conditional results were obtained by Murty--Pasten \cite{MP2} and by Pasten \cite{Pasten}. We establish the same result unconditionally.

\begin{theorem}
\label{tMain1}
Let $R$ be a finitely generated ring over $\Z$ with $|R| = \infty$. Then Hilbert's tenth problem has a negative answer for $R$.
\end{theorem}

We prove Theorem \ref{tMain1} by proving the following conjecture of Denef and Lipshitz \cite{DL}.

\begin{theorem}
\label{tMain2}
Let $K$ be a number field. Then $\Z$ is Diophantine over $O_K$.
\end{theorem}

Eisentr\"ager \cite{Eis} proved in her PhD thesis that Theorem \ref{tMain2} implies Theorem \ref{tMain1}, so it remains to establish Theorem \ref{tMain2}. 

Approximately two months after the release of our initial preprint, Alp\"oge--Bhargava--Ho--Shnidman \cite{ABHS} proved a rank growth result of similar strength as our Theorem \ref{tRelative}. Their main result recovers our Theorem \ref{tMain1} and Theorem \ref{tMain2}.

\subsection{Additive combinatorics}
Our main new insight is to bring additive combinatorics into this field. For now, let $L/K$ be a quadratic extension. Our goal will be to construct an elliptic curve $E$ over $K$ with
\begin{align}
\label{eRankGrowth}
\mathrm{rk} \, E(L) = \mathrm{rk} \, E(K) > 0.
\end{align}
In fact, we may and will later assume quite a bit more about $L/K$ to aid with some of the technical arguments, see Theorem \ref{tRelative} for a precise statement. 

Mazur--Rubin \cite{MR} conditionally construct elliptic curves satisfying \eqref{eRankGrowth} by taking an elliptic curve $E$ and quadratic twisting it to control the $2$-Selmer rank. In their work they take an elliptic curve with $\Gal(K(E[2])/K) \cong S_3$ as large as possible. Then they rely on finiteness of Sha to go back to ranks.

We instead perform $2$-descent on elliptic curves with full rational $2$-torsion. Let $E$ be the elliptic curve
$$
y^2 = (x - a_1) (x - a_2) (x - a_3).
$$
Our first idea is to consider the twist by $t := d(c - a_1d)(c - a_2d)(c - a_3d)$, so $E^t$ has the equation
$$
ty^2 = (x - a_1) (x - a_2) (x - a_3),
$$
which tautologically has the rational point $(x, y) = (c/d, 1/d^2)$. This point is readily seen to be non-torsion for all but finitely many twists $t$ and hence the rank of $E^t$ is positive. Thus, if we were able to show that the $2$-Selmer rank of $E^t$ over $L$ equals $2 + \mathrm{rk} \, E(K)$, then we would have unconditionally established equation \eqref{eRankGrowth}. However, as the Selmer growth results in Morgan--Paterson \cite{MP} suggest, this is typically not possible so we have opted to prove instead that the $2$-Selmer rank of $E^{\alpha t}$ is $2$, where $\alpha$ is such that $L = K(\sqrt{\alpha})$.

Our second idea is to take advantage of additive combinatorics. Indeed, in general it is exceptionally difficult to control Selmer ranks of polynomial twists; there are no distributional results known outside linear polynomials. However, the change of Selmer ranks is well understood for prime elements. Hence we consider the situation that $d$, $c - a_1d$, $c - a_2d$ and $c - a_3d$ are all prime elements. This is possible by recent work of Kai \cite{Kai}, which is a Green--Tao type result for number fields based on the earlier works \cite{GT1, GT2, GT3, GTZ}. In fact, we only need an inverse theorem for the Gowers $U^3$ norm, which is technically simpler and established previously in \cite{GT0}. 

\subsection{Layout of the paper}
In Section \ref{sRed} we perform a standard reduction step to reduce the problem to establishing equation \eqref{eRankGrowth} for certain quadratic extensions $L/K$. This reduction step allows one to assume that $K$ contains any given totally real field. We will use this flexibility merely to ensure that $K$ itself has at least $32$ real places, which we will take advantage of in Section \ref{sGT}. Moreover, unlike Mazur--Rubin \cite{MR} we reduce to the case where $L/K$ is quadratic instead of cyclic. That this suffices was first noticed by Shlapentokh, see \cite[Theorem 4.8]{MRS}.

The main result of Section \ref{sRed} is that Theorem \ref{tRelative} implies Theorem \ref{tMain2}. The remainder of the paper is dedicated to establishing Theorem \ref{tRelative}.

In Section \ref{sBackground} we develop some background theory. We perform $2$-descent over general number fields, and give several results to compute the parity of $2$-Selmer. Moreover, we construct suitable ``starting'' elliptic curves that we shall twist in the later sections. We end this section by giving a reasonable parametrization of $H^1(G_K, \mathbb{F}_2)$.

In Section \ref{sGen} we define a transition process that computes the change in $2$-Selmer rank of $E^\kappa$ to $E^{\kappa \pi}$. Since this transition is rather delicate, we spend a great amount of time to carefully control the $2$-Selmer rank of $E^\kappa$ in Section \ref{sGT}. 

After having suitably prepared our Selmer groups for our additive combinatorics arguments, we finish the argument in Section \ref{sGT}. At this stage a final obstacle arises as one can not independently prescribe $c - a_1d$, $c - a_2d$, $c - a_3d$ modulo $d$. To work around this issue, we make use of infinite places to carefully control some terms arising from quadratic reciprocity.

\subsection*{Acknowledgements}
This project was initiated when the second author visited the first author at the Institute for Theoretical Studies (ETH Z\"urich), and both authors wish to thank the institute for its financial support and excellent working conditions. The first author gratefully acknowledges the support of Dr. Max R\"ossler, the Walter Haefner Foundation and the ETH Z\"urich Foundation. The first author also acknowledges the support of the Dutch Research Council (NWO) through the Veni grant ``New methods in arithmetic statistics''.

The manuscript was completed while the authors were at the workshop \emph{Nilpotent counting problems in arithmetic statistics} at AIM. We are very grateful to AIM and the organizers Brandon Alberts, Yuan Liu and Melanie Wood.

The authors are deeply indebted to Efthymios Sofos, who answered many questions about additive combinatorics during the project. Without his constant support, the project would most likely not have come to be. We would also like to thank Gunther Cornelissen, Christopher Frei, Andrew Granville, Ben Green, Wataru Kai, Emmanuel Kowalski, Vivian Kuperberg, Jeff Lagarias and Alexander Smith for helpful conversations. We are also indebted to Hector Pasten for explaining us how to reduce Hilbert's tenth problem to rank growth in quadratic extensions instead of general cyclic extensions.

We would like to thank Adam Morgan for valuable feedback on Section \ref{sBackground}. Moreover, we are most grateful to him for the shared time at the Max Planck Institute in Bonn, where he explained to us how to control the $2$-Selmer group with the transition process from Section \ref{sGen}. In particular, his excellent notes formed the basis of our Section \ref{sGen}.

Finally, we would like to express our admiration for the dedication of the anonymous referees. Their careful reading and their numerous suggestions have greatly improved this manuscript. In particular, their suggestions have led to several simplifications of the original arguments.
\section{Initial reductions}
\label{sRed}
In this section we reduce our main theorem to establishing a certain rank growth condition of elliptic curves. This reduction is standard in the literature, see for example \cite[Section 8]{MR}. We redo the arguments here as we are going to impose some extra conditions on the relevant extensions of number fields that will ease some future arguments. We start by recalling three well-known results.

\begin{lemma}
\label{lMR}
Suppose $K \subseteq L$ are number fields. Then:
\begin{itemize}
\item[(i)] If $D_1, D_2 \subseteq O_L$ are Diophantine over $O_L$, then so is $D_1 \cap D_2$.
\item[(ii)] If $D \subseteq O_K$ is Diophantine over $O_K$, and $O_K$ is Diophantine over $O_L$, then $D$ is Diophantine over $O_L$.
\item[(iii)] If $\Z$ is Diophantine over $O_L$, then $\Z$ is Diophantine over $O_K$.
\end{itemize}
\end{lemma}

\begin{proof}
See \cite[Lemma 8.3]{MR}.
\end{proof}

\begin{lemma}
\label{lReal}
If $M$ is a totally real number field, then $\Z$ is Diophantine over $O_M$.
\end{lemma}

\begin{proof}
This is the statement of \cite[Theorem]{Denef}.
\end{proof}

\begin{lemma}
\label{lPoonen}
If $L/K$ is an extension of number fields and if there exists an elliptic curve $E/K$ such that $\mathrm{rk} \, E(K) = \mathrm{rk} \, E(L) > 0$, then $O_K$ is Diophantine over $O_L$.
\end{lemma}

\begin{proof}
See \cite[Theorem 1.9]{Shl}.
\end{proof}

Having stated these results, we will now perform our reduction step. Let us first state the theorem that we shall prove in the remaining sections of the paper.

\begin{theorem}
\label{tRelative}
Let $K$ be a number field and set $L := K(i)$. Assume that $K$ has at least $32$ real places. Then there exists an elliptic curve $E$ over $K$ such that $\mathrm{rk} \, E(K) = \mathrm{rk} \, E(L) > 0$.
\end{theorem}

We shall prove Theorem \ref{tRelative} in Section \ref{sGT}. In order to prove that Theorem \ref{tRelative} implies Theorem \ref{tMain2}, we shall pass through the following intermediate result.

\begin{corollary}
\label{cRelative2}
Let $K$ be a number field and set $L := K(i)$. Assume that 
\begin{itemize}
\item[(i)] $L$ is Galois over $\Q$,
\item[(ii)] $K$ has a real place,
\item[(iii)] $K$ contains $K_0 := \Q(\sqrt{5}, \sqrt{7}, \sqrt{11}, \sqrt{13}, \sqrt{17}, \sqrt{19})$.
\end{itemize}
Then there exists an elliptic curve $E$ over $K$ such that $\mathrm{rk} \, E(K) = \mathrm{rk} \, E(L) > 0$.
\end{corollary}

\begin{proof}[Proof that Theorem \ref{tRelative} implies Corollary \ref{cRelative2}]
Let $K$ be a number field, set $L := K(i)$ and assume that conditions $(i), (ii), (iii)$ hold.

Write $\mathrm{C}_G(S)$ for the centralizer of a subset $S$ inside a group $G$. Fix an embedding $\sigma: K \rightarrow \mathbb{R}$ and extend it to $\sigma_1: L \rightarrow \mathbb{C}$. We will start by proving that 
\begin{equation}
\label{eTauClaim}
\{\tau \in \Gal(L/\Q) : \sigma_1(\tau(K)) \subseteq \mathbb{R}\} = \mathrm{C}_{\Gal(L/\Q)}(\Gal(L/K)). 
\end{equation}
Write $\langle \rho \rangle = \Gal(L/K)$. Denoting complex conjugation by $z \mapsto \bar{z}$, we observe that $\sigma_1 \circ \rho = \overline{\sigma_1}$. Hence if $\tau \in \mathrm{C}_{\Gal(L/\Q)}(\langle \rho \rangle)$, we have
$$
\overline{\sigma_1(\tau(\alpha))} = \sigma_1(\rho(\tau(\alpha))) = \sigma_1(\tau(\rho(\alpha))) = \sigma_1(\tau(\alpha))
$$
for all $\alpha \in K$, so $\sigma_1(\tau(K)) \subseteq \mathbb{R}$. Conversely, if $\sigma_1(\tau(K)) \subseteq \mathbb{R}$, then
$$
\sigma_1(\rho(\tau(\alpha))) = \overline{\sigma_1(\tau(\alpha))} = \sigma_1(\tau(\alpha)) = \sigma_1(\tau(\rho(\alpha)))
$$
for all $\alpha \in K$, hence $[\rho, \tau]$ must fix $K$. This implies $[\rho, \tau] \in \Gal(L/K)$, therefore $[\rho, \tau]$ is the identity (as the alternative $[\rho, \tau] = \rho$ forces $\rho$ to be the identity, contradiction). This proves equation \eqref{eTauClaim}.

In order to apply Theorem \ref{tRelative}, we need to show that $K$ has at least $32$ real places. Note that there is a $2$-to-$1$ surjective map from $\{\tau \in \Gal(L/\Q) : \sigma_1(\tau(K)) \subseteq \mathbb{R}\}$ to the set of real embeddings of $K$. Therefore, by equation \eqref{eTauClaim}, the number of real embeddings of $K$ is exactly $|\mathrm{C}_{\Gal(L/\Q)}(\langle \rho \rangle)|/2$. Let $C$ be the conjugacy class of $\rho$ in $\Gal(L/\Q)$. Then $C$ is contained in the normal subgroup $\Gal(L/K_0)$, so $|C| \leq [L : K_0]$. By the orbit–stabilizer formula, $|\mathrm{C}_{\Gal(L/\Q)}(\langle \rho \rangle)| = |\Gal(L/\Q)|/|C| \geq [L : \Q]/[L : K_0] = [K_0 : \Q] = 64$.
\end{proof}

We shall now prove that Corollary \ref{cRelative2} implies Theorem \ref{tMain2} (stated here again as Theorem \ref{tHilbert10}). The remainder of this paper is dedicated to establishing Theorem \ref{tRelative}.

\begin{theorem}
\label{tHilbert10}
Let $F$ be a number field. Then $\Z$ is Diophantine over $O_F$.
\end{theorem}

\begin{proof}[Proof that Corollary \ref{cRelative2} implies Theorem \ref{tHilbert10}]
Let $F$ be any number field. Define $L$ to be the compositum of the normal closure of $F$ together with $K_0(i)$. Let $K$ be any field such that $L = K(i)$ and such that $K$ has at least one real place. Then $K$ must contain all totally real subfields of $L$, hence contains $K_0$. We conclude that $O_K$ is Diophantine over $O_L$ by Lemma \ref{lPoonen} and Corollary \ref{cRelative2} applied to $L/K$. 

Take $M$ to be the intersection of all $K$ that satisfy $L = K(i)$ and have at least one real place. As we have just proven that $O_K$ is Diophantine over $O_L$ for all such $K$, then so is $O_M$ by Lemma \ref{lMR}$(i)$. Note that $M$ is Galois over $\Q$, and since there exists a field $K$ containing $M$ with at least one real place, $M$ must be totally real. So $\Z$ is Diophantine over $O_M$ by Lemma \ref{lReal}, and thus $\Z$ is Diophantine over $O_L$ by Lemma \ref{lMR}$(ii)$. This implies that $\Z$ is Diophantine over $O_F$ by Lemma \ref{lMR}$(iii)$.
\end{proof}
\section{Background theory}
\label{sBackground}

\subsection{Elliptic curves}
Given a field $K$ and $a_1, a_2, a_3 \in K$, we define 
\begin{equation}
\label{eabc}
\alpha := a_1 - a_2, \quad \quad \beta := a_1 - a_3, \quad \quad \gamma := a_2 - a_3. 
\end{equation}
Write $(-, -)_{\mathrm{Weil}} : E[2] \times E[2] \rightarrow \mu_2$ for the Weil pairing.

\begin{lemma} 
\label{lemma: delta map formula}
Let $K$ be a field of characteristic $0$. Let $a_1, a_2, a_3 \in K$ be distinct and let $E$ be the elliptic curve defined by
$$
E: y^2 = (x - a_1) (x - a_2) (x - a_3).
$$
Let $P_1 = (a_1, 0)$ and $P_2 = (a_2, 0)$. Let $\lambda_i$ be the homomorphism $E[2] \rightarrow \mu_2$ given by $P \mapsto (P, P_i)_{\mathrm{Weil}}$. Under the isomorphism $(\lambda_1, \lambda_2) : H^1(G_K, E[2]) \cong (K^\ast/K^{\ast 2})^2$, the Kummer map $\delta: E(K)/2E(K) \xhookrightarrow{} H^1(G_K, E[2]) \cong (K^\ast/K^{\ast 2})^2$ is given explicitly by
$$
P = (x, y) \mapsto
\begin{cases}
(x - a_1, x - a_2) &\textup{if } P \neq P_1, P_2 \\
(\alpha \beta, \alpha) &\textup{if } P = P_1 \\
(-\alpha, -\alpha \gamma) &\textup{if } P = P_2.
\end{cases}
$$
\end{lemma}

\begin{proof}
See \cite[Proposition 1.4, Chapter 10]{Silverman}.
\end{proof}

Throughout the paper we shall often implicitly identify $H^1(G_K, E[2])$ with $(K^\ast/K^{\ast 2})^2 \cong H^1(G_K, \mathbb{F}_2)^2$. Whenever we do so, the implicit identification used is $(\lambda_1, \lambda_2)$.

In order to show that the rational point that we will construct later is not torsion, the following fact will be useful.

\begin{lemma}
\label{lTorsion}
Let $K$ be a number field and let $E$ be an elliptic curve defined over $K$. For all but finitely many $d \in K^\ast/K^{\ast 2}$, we have
$$
E^d(K)^{\textup{tors}} = E^d(K)[2].
$$
\end{lemma}

\begin{proof}
Recall that for each $d \in K^\ast$ the map 
$$
(x, y) \mapsto (x, \sqrt{d} \cdot y)
$$
induces an injective group homomorphism $\phi_d: E^d(K) \to E(\overline{K})$. The image of $\phi_d$ consists of the zero element together with those points $(x, y) \in E(\overline{K})$ satisfying $x \in K$ and $y \in K \sqrt{d}$. Therefore, if we take $d_1, d_2 \in K^\ast$ with different classes in $K^\ast/K^{\ast 2}$, we have that
$$
\text{im}(\phi_{d_1}) \cap \text{im}(\phi_{d_2}) = E(K)[2].
$$
Hence it suffices to show that
$$
\left| \bigcup_{d \in K^\ast/K^{\ast 2}} \text{im}(\phi_d)^{\text{tors}} \right| < \infty.
$$
However, $\bigcup_{d \in K^\ast/K^{\ast 2}} \text{im}(\phi_d)^{\text{tors}} \subseteq E(\overline{K})^{\text{tors}}$ is a subset of bounded height and degree, thus we conclude with Northcott's theorem. 
\end{proof}

%We claim that for any fixed odd prime $p$ there are at most $2$ values of $d \in K^\ast/K^{\ast 2}$ such that
%$$
%E^d(K)[p] \neq 0.
%$$
%Indeed, assume for the sake of contradiction that there are distinct elements $d_1, \dots, d_3 \in K^\ast/K^{\ast 2}$ with $E^{d_i}(K)[p] \neq 0$. We write $E' = E^{d_1}$. Since $|\Gal(K(\sqrt{d_1 d_2}, \sqrt{d_1 d_3})/K)| = 4$ is coprime with $p$, we have the decomposition
%$$
%E'(K(\sqrt{d_1 d_2}, \sqrt{d_1 d_3}))[p] = E^{d_1}(K)[p] \oplus E^{d_2}(K)[p] \oplus E^{d_3}(K)[p] \oplus E^{d_1d_2d_3}(K)[p].
%$$
%But such a decomposition is clearly impossible, as the LHS has dimension at most $2$, while the RHS has dimension at least $3$.

%By Merel's theorem, there is a finite subset of the primes $\mathcal{P}$ depending only on the degree $[K : \Q]$ of $K$ such that
%$$
%E^d(K)[p] = 0 \quad \quad \text{ for all } p \not \in \mathcal{P}.
%$$
%Excluding the at most $2|\mathcal{P}|$ twists $d$ with $E^d(K)[p] \neq 0$ for some odd prime $p \in \mathcal{P}$ and excluding all $d$ with $K(\sqrt{d})/K$ ramified only at $2$, infinite places and the bad places of $E$ finishes the proof thanks to Lemma \ref{lTorsion2}.

\begin{lemma}
\label{lRankDecomposition}
Let $K$ be a number field, let $d \in K^\ast/K^{\ast 2}$ be non-trivial and set $L := K(\sqrt{d})$. Let $E$ be an elliptic curve defined over $K$. Then
$$
\textup{rk} \, E(L) = \textup{rk} \, E(K) + \textup{rk} \, E^d(K).
$$
\end{lemma}

\begin{proof}
This is a standard result, which follows from decomposing $E(L) \otimes_{\Z} \Q$ into its plus and minus part under the action of $\Gal(L/K)$.
\end{proof} 

\subsection{Root numbers}
The global root number $w(E/K)$ is defined by
$$
w(E/K) = \prod_v w(E/K_v),
$$
where $w(E/K_v)$ is the local root number. We shall not need the precise definition of this local root number as our next theorem provides a convenient way to calculate it.

\begin{theorem}[Root number facts]
\label{tRoot}
Let $E$ be an elliptic curve defined over a local field $\mathcal{K}$ of characteristic $0$. Let $k$ be the residue field of $\mathcal{K}$ and let $v$ be the associated normalized valuation of $\mathcal{K}$. Then we have:
\begin{itemize}
\item[(i)] If $E/\mathcal{K}$ has good reduction, then $w(E/\mathcal{K}) = +1$.
\item[(ii)] If $E/\mathcal{K}$ has additive, potentially good reduction and $\mathrm{char}(k) \geq 5$, then 
$$
w(E/\mathcal{K}) = (-1)^{\lfloor \frac{v(\Delta) |k|}{12} \rfloor}.
$$
\item[(iii)] If $\mathcal{K}$ is archimedean, then $w(E/\mathcal{K}) = -1$.
\end{itemize}
\end{theorem}

\begin{proof}
See \cite[Theorem 2.3]{CD} or \cite[Theorem 3.1]{DD2}. These statements were originally proven in Kobayashi \cite{Kobayashi} and Rohrlich \cite{Rohrlich}.
\end{proof}

The relevance of these root number computations is that we can relate them to the parity of the $2$-Selmer rank by our next theorem. Recall that $\mathrm{Sel}^{2^{\infty}}(E/K) := \varinjlim_n \mathrm{Sel}^{2^n}(E/K)$. Define $\mathrm{rk}_2(E/K)$ to be the $\Z_2$-corank of the $2^\infty$-Selmer group.

\begin{theorem}[2-parity]
\label{tCes}
Let $K$ be a number field and let $E/K$ be an elliptic curve with $E(K)[2] \cong \mathbb{F}_2^2$. We have
$$
(-1)^{\mathrm{rk}_2(E/K)} = w(E/K).
$$
\end{theorem}

\begin{proof}
This was originally proven in \cite[Theorem 1.8]{DD} and generalized in \cite[Theorem 1.4]{Ces}.
\end{proof}

% We will make use of the following auxiliary fact. 

% \begin{proposition} 
% \label{prop: groups with alternating pairings are squares}
% Let $A$ be a finite abelian group equipped with a perfect non-degenerate alternating pairing $b:A \times A \to \frac{\Q}{\Z}$. Then there exists a finite abelian group $B$ such that
% $$
% A \simeq_{\textup{ab.gr.}} B \oplus B.
% $$
% \end{proposition}

% \begin{proof}
% See \cite{Davydov}[Lemma 5.2].   
% \end{proof}

\begin{corollary}
\label{cCes}
Let $K$ be a number field and let $E/K$ be an elliptic curve. Suppose that $E(K)[2] \cong \mathbb{F}_2^2$. Then
$$
(-1)^{\dim_{\FF_2} \mathrm{Sel}^2(E/K)} = w(E/K).
$$
\end{corollary}

\begin{proof}
In virtue of Theorem \ref{tCes}, it suffices to prove that $\dim_{\mathbb{F}_2} \mathrm{Sel}^2(E/K)$ and $\mathrm{rk}_2(E/K)$ are congruent modulo $2$, which we do below. We start with the exact sequences
\begin{gather*}
0 \to E(K)/2E(K) \to \mathrm{Sel}^2(E/K) \to \Sha(E/K)[2] \to 0 \\
0 \to E(K) \otimes_{\Z} \frac{\Q_2}{\Z_2} \to \mathrm{Sel}^{2^{\infty}}(E/K) \to \Sha(E/K)[2^\infty] \to 0.
\end{gather*}
Write $\Sha(E/K)_{\text{div}}$ for the largest divisible subgroup of $\Sha(E/K)$. The second exact sequence implies that $\mathrm{rk}_2(E/K) = \mathrm{rk} \, E(K) + \dim_{\mathbb{F}_2} \Sha(E/K)_{\text{div}}[2]$. Recalling that $E(K)[2] \cong \mathbb{F}_2^2$ by assumption and taking dimensions in the first exact sequence above gives
\begin{align*}
\dim_{\mathbb{F}_2} \mathrm{Sel}^2(E/K) &= \dim_{\mathbb{F}_2} E(K)[2] + \mathrm{rk} \, E(K) + \dim_{\mathbb{F}_2} \Sha(E/K)[2] \\
&= 2 + \mathrm{rk} \, E(K) + \dim_{\mathbb{F}_2} \Sha(E/K)_{\text{div}}[2] + \dim_{\mathbb{F}_2} \frac{\Sha(E/K)}{\Sha(E/K)_{\text{div}}}[2] \\
&= 2 + \mathrm{rk}_2(E/K) + \dim_{\mathbb{F}_2} \frac{\Sha(E/K)}{\Sha(E/K)_{\text{div}}}[2].
\end{align*}
By \cite[Corollary 12]{PS}, the Poonen--Stoll class pairs trivially with itself. Then \cite[Theorem 8, 1 implies 3]{PS} shows that $\dim_{\mathbb{F}_2} \frac{\Sha(E/K)}{\Sha(E/K)_{\text{div}}}[2]$ is even, which immediately gives the corollary.
\end{proof}

Finally, we are ready to construct a starting curve that we will twist later. Since our methods are fairly flexible, we shall not need to demand too much on our starting curve. Given three distinct real numbers $x_1, x_2, x_3$, we define $p(x_1, x_2, x_3)$ to be the unique permutation $\sigma \in S_3$ such that
$$
x_{\sigma(1)} < x_{\sigma(2)} < x_{\sigma(3)}.
$$
Given $a_1, a_2, a_3$, we recall the notation $\alpha = a_1 - a_2$, $\beta = a_1 - a_3$ and $\gamma = a_2 - a_3$.

\begin{lemma}
\label{lStartingCurve}
Let $K$ be a number field with at least $25$ real embeddings. Then there exist elements $a_1, a_2, a_3 \in O_K$ such that 
\begin{enumerate}
\item[(a)] $-1, \alpha, \beta, \gamma \in K^\ast$ are linearly independent in $K^\ast/K^{\ast 2}$;
\item[(b)] there exist real embeddings $\sigma_1, \dots, \sigma_{24}$ such that for every $f \in S_3$ there are precisely four indices $i$ satisfying 
$$
p(\sigma_i(a_1), \sigma_i(a_2), \sigma_i(a_3)) = f;
$$
\item[(c)] the elliptic curve 
$$
E: y^2 = (x - a_1) (x - a_2) (x - a_3)
$$
satisfies $w(E/K) = +1$.
\end{enumerate}
\end{lemma}

\begin{proof}
Fix some real embeddings $\sigma_1, \dots, \sigma_{25}$. There exist $a_1', a_2', a_3' \in O_K$ satisfying $(b)$ as $O_K$ is naturally a lattice of full rank inside $K \otimes_{\Z} \mathbb{R}$ and thus meets every orthant.

Let $\widetilde{E}$ be the elliptic curve given by
$$
\widetilde{E}: y^2 = (x - a_1') (x - a_2') (x - a_3').
$$
Write $N$ for the product of finite bad places of $\widetilde{E}$. By the Mitsui prime ideal theorem (see \cite[Theorem, p.~35]{Mitsui} or \cite[Theorem 1.2.1]{KaiM} for a more modern version), there exists an odd principal prime element $\pi$ of good reduction that has the same sign as $w(\widetilde{E}/K)$ at $\sigma_{25}$, is positive at all other real places and satisfies $\pi \equiv 1 \bmod 24N$. Applying Hilbert reciprocity to the symbol $(\pi, -1)$, we see that $|O_K/\pi| \equiv 1 \bmod 4$ if $w(\widetilde{E}/K) = +1$ and $|O_K/\pi| \equiv 3 \bmod 4$ if $w(\widetilde{E}/K) = -1$. We now take $E := \widetilde{E}^\pi$, $a_1 := \pi a_1'$, $a_2 := \pi a_2'$ and $a_3 := \pi a_3'$. 

Because $\widetilde{E}$ has good reduction at the prime $(\pi)$, the $\pi$-adic valuation of $\Delta_{\widetilde{E}}$ equals zero by \cite[Proposition 5.1, Chapter 7]{Silverman}. Then a direct inspection of the formula for $\Delta_E$ shows that the $\pi$-adic valuation of $\Delta_E$ equals $6$. In particular, Theorem \ref{tRoot} gives that $w(E/K_{(\pi)}) = -1$. Since $\widetilde{E}$ and $E$ are isomorphic at all finite places of bad reduction for $\widetilde{E}$ and are also isomorphic at all $2$-adic and $3$-adic places, it follows from Theorem \ref{tRoot} that $w(E/K) = +1$.

Since $a_1', a_2', a_3'$ satisfy $(b)$ and since $\pi$ is totally positive at $\sigma_1, \dots, \sigma_{24}$, it follows that $a_1, a_2, a_3$ also satisfy $(b)$. Finally, $(a)$ is in fact a consequence of $(b)$.
\end{proof}
\section{Selmer ranks as transition processes}
\label{sGen}
Let $K$ be a number field and let $E/K$ be an elliptic curve with $E(K)[2] \cong \mathbb{F}_2^2$. The central result of this section is a description of $2$-Selmer ranks of quadratic twists of $E$ in terms of a transition process. We denote by $\Omega_K$ the set of places of $K$. Take $T$ to be any finite subset of $\Omega_K$ containing all $2$-adic places of $K$, all infinite places of $K$, all places of bad reduction of $E$ and all places where $\alpha$, $\beta$ or $\gamma$ has non-zero valuation; for the definition of $\alpha$, $\beta$ and $\gamma$, see equation \eqref{eabc}. 

To describe the Selmer rank as a transition process, we will introduce a sequence of Selmer structures. Let $A$ be a finite $G_K$-module. For each $v \in \Omega_K$, we fix an embedding $G_{K_v} \hookrightarrow G_K$ and we write $I_v$ for the inertia subgroup of $G_{K_v}$. The embedding $G_{K_v} \hookrightarrow G_K$ gives rise to a natural restriction map $\res_v: H^1(G_K, A) \rightarrow H^1(G_{K_v}, A)$. Define
$$
H^1_{\text{nr}}(G_{K_v}, A) := \ker\left(H^1(G_{K_v}, A) \rightarrow H^1(I_v, A)\right).
$$
A Selmer structure $\mathcal{L} = (\mathcal{L}_v)_{v \in \Omega_K}$ (in the language of \cite[Definition 8.7.8]{NSW}, this is called a collection of local conditions) is a collection of subgroups $\mathcal{L}_v \subseteq H^1(G_{K_v}, A)$ such that $\mathcal{L}_v = H^1_{\text{nr}}(G_{K_v}, A)$ for all but finitely many $v$. To a Selmer structure $\mathcal{L}$, we associate its Selmer group via the formula
$$
\mathrm{Sel}_{\mathcal{L}}(G_K, A) := \ker\left(H^1(G_K, A) \xrightarrow{\prod_v \res_v} \prod_{v \in \Omega_K} \frac{H^1(G_{K_v}, A)}{\mathcal{L}_v}\right).
$$
Let $n \in \Z_{\geq 0}$ and let $v_1, \dots, v_n$ be distinct places outside of $T$. For each $1 \leq j \leq n$, let $\pi_j$ be an element of $K_{v_j}^\ast/K_{v_j}^{\ast 2}$ with odd valuation (note that there are exactly two such choices) and write $\boldsymbol{\pi} = (\pi_j)_{1 \leq j \leq n}$ for the resulting vector. For each integer $0 \leq i \leq n$, we define a Selmer structure $\mathcal{L}_{i, \boldsymbol{\pi}} = (\mathcal{L}_{i, \boldsymbol{\pi}, v})_{v \in \Omega_K}$ by
\begin{align}
\label{eSelmerStructure}
\mathcal{L}_{i, \boldsymbol{\pi}, v} = 
\begin{cases}
\delta(E(K_v)) &\text{if } v \in T \\
\langle (\alpha \beta, \pi_j \alpha), (-\pi_j \alpha, -\alpha \gamma) \rangle &\text{if } v = v_j \text{ for some } j \in \{1, \dots, i\} \\
H^1_{\text{nr}}(G_{K_v}, E[2]) &\text{otherwise.}
\end{cases}
\end{align}
We shall see in Lemma \ref{lEasy} that $\langle (\alpha \beta, \pi_j \alpha), (-\pi_j \alpha, -\alpha \gamma) \rangle = \delta(E^{\pi_j}(K_v))$, which is the main motivation for defining the second line of equation \eqref{eSelmerStructure} as above. 

For $0 \leq i \leq n - 1$, define
$$
n_i := \dim_{\mathbb{F}_2} \mathrm{Sel}_{\mathcal{L}_{i + 1, \boldsymbol{\pi}}}(G_K, E[2]) - \dim_{\mathbb{F}_2} \mathrm{Sel}_{\mathcal{L}_{i, \boldsymbol{\pi}}}(G_K, E[2]).
$$
In order to study $n_i$, we introduce a second Selmer structure $\mathcal{L}_{i, \boldsymbol{\pi}}' = (\mathcal{L}_{i, \boldsymbol{\pi}, v}')_{v \in \Omega_K}$ through
$$
\mathcal{L}_{i, \boldsymbol{\pi}, v}' = 
\begin{cases}
\delta(E(K_v)) &\text{if } v \in T \\
\langle (\alpha \beta, \pi_j \alpha), (-\pi_j \alpha, -\alpha \gamma) \rangle &\text{if } v = v_j \text{ for some } j \in \{1, \dots, i\} \\
H^1(G_{K_v}, E[2]) &\text{if } v = v_{i + 1} \\
H^1_{\text{nr}}(G_{K_v}, E[2]) &\text{otherwise}
\end{cases}
$$
and
\begin{align}
\label{eRelaxedA}
A_{v_{i + 1}} := \res_{v_{i + 1}} \mathrm{Sel}_{\mathcal{L}_{i, \boldsymbol{\pi}}'}(G_K, E[2]).
\end{align}
As we shall now summarize, an important insight of Poonen--Rains \cite{PR} is that the local Tate pairing naturally comes from a certain quadratic form $q_{E, v}: H^1(G_{K_v}, E[2]) \rightarrow \mathrm{Br}(K_v)[2]$, which was first constructed by Zarhin \cite{Zarhin}. We will then use this theory to show that $A_{v_{i + 1}}$ is a maximal isotropic subspace (see Proposition \ref{prop: relaxed Selmer image maximal isotropic}).

We start by recalling some basic properties of this construction. Let $v \in \Omega_K$. Identifying $E[2]$ with its dual $G_{K_v}$-module through the Weil pairing, we see that local Tate duality induces an alternating non-degenerate bilinear pairing
$$
b_{E, v}: H^1(G_{K_v}, E[2]) \times H^1(G_{K_v}, E[2]) \to \text{Br}(K_v)[2].
$$

\begin{proposition} 
\label{prop: loc.glob.quad.form}
$(a)$ For each place $v$ of $K$, there exists a quadratic form 
$$
q_{E, v}: H^1(G_{K_v}, E[2]) \to \textup{Br}(K_v)[2]
$$ 
inducing the pairing $b_{E, v}$. Furthermore, $\delta(E(K_v)/2E(K_v))$ is a maximal isotropic subspace of the quadratic space $(H^1(G_{K_v}, E[2]), q_{E, v})$.

$(b)$ There exists a quadratic form $q_{E, K}: H^1(G_K, E[2]) \to \textup{Br}(K)[2]$ satisfying 
$$
\res_v \circ q_{E, K} = q_{E, v} \circ \res_v.
$$
\end{proposition}

\begin{proof}
$(a)$ The construction of the quadratic forms $q_{E,v}$ and $q_{E,K}$ can be found in \cite[Section 4.1]{PR}. The fact that it induces $b_{E, v}$ is given in \cite[Corollary 4.7]{PR}, while the last part is \cite[Proposition 4.11]{PR}. 

$(b)$ The forms $q_{E, K}$ and $q_{E, v}$ in \cite[Section 4.1]{PR} are constructed as connecting homomorphisms.
% of a sequence of respectively $G_K$-modules and $G_{K_v}$-modules. 
Hence the desired claim follows from the fact that connecting homomorphisms commute with restriction to subgroups.
\end{proof}

We will now use Proposition \ref{prop: loc.glob.quad.form} to prove several important equalities between the local conditions appearing in the definitions of $\mathcal{L}_{i, \boldsymbol{\pi}}$ and $\mathcal{L}_{i, \boldsymbol{\pi}}'$.

\begin{lemma}
\label{lEasy}
\begin{enumerate}
\item[(a)] If $v \in \Omega_K$ is odd, then
$$
\dim_{\FF_2} H^1(G_{K_v}, E[2]) = \dim_{\FF_2} H^1(G_{K_v}, \mathbb{F}_2^2) = 2\dim_{\FF_2} H^1(G_{K_v}, \mathbb{F}_2) = 4.
$$
Moreover, we have 
$$
\dim_{\FF_2} H^1_{\textup{nr}}(G_{K_v}, E[2]) = 2 \dim_{\FF_2} H^1_{\textup{nr}}(G_{K_v}, \mathbb{F}_2) = 2.
$$
\item[(b)] If $v \in \Omega_K$ is odd and of good reduction for $E$, then
$$
\delta(E(K_v)/2E(K_v)) = H^1_{\textup{nr}}(G_{K_v}, E[2]).
$$
\item[(c)] If $v \in \Omega_K - T$ and if $\pi$ is a uniformizer for $K_v$, then
$$
\delta(E^\pi(K_v)/2E^\pi(K_v)) = \langle (\alpha \beta, \pi \alpha), (-\pi \alpha, -\alpha \gamma) \rangle.
$$
\item[(d)] If $v \in \Omega_K - T$ and if $\pi$ is a uniformizer for $K_v$, then
$$
\langle (\alpha \beta, \pi \alpha), (-\pi \alpha, -\alpha \gamma) \rangle \cap H^1_{\textup{nr}}(G_{K_v}, E[2]) = 0.
$$
\end{enumerate}
\end{lemma}

\begin{proof}
The first half of part $(a)$ is immediate from the explicit description of $K_v^\ast/K_v^{\ast 2}$ for local fields $K_v$. The second half of part $(a)$ follows from the isomorphism $E[2] \cong \mathbb{F}_2^2$ as Galois modules and the fact $\dim_{\FF_2} H^1_{\textup{nr}}(G_{K_v}, \mathbb{F}_2) = 1$. Part $(b)$ follows from \cite[Remark 4.12 and Proposition 4.13]{PR}.

As for part $(c)$, we start by observing that $\delta(E^\pi(K_v)/2E^\pi(K_v))$ is a maximal isotropic subspace of $q_{E^\pi, v}$ by Proposition \ref{prop: loc.glob.quad.form}$(a)$. Therefore it follows from part $(a)$ that 
$$
\dim_{\FF_2} \delta(E^\pi(K_v)/2E^\pi(K_v)) = 2. 
$$
Inspecting valuations and using that $v \not \in T$, we see that $\dim_{\FF_2} \langle (\alpha \beta, \pi \alpha), (-\pi \alpha, -\alpha \gamma) \rangle = 2$. Since $(\alpha \beta, \pi \alpha)$ and $(-\pi \alpha, -\alpha \gamma)$ are inside $\delta(E^\pi(K_v)/2E^\pi(K_v))$ by Lemma \ref{lemma: delta map formula}, part $(c)$ follows. Part $(d)$ is a direct consequence of inspecting valuations.
\end{proof}

We now describe the initial state of the transition process, which is another important motivation for the definition of $\mathcal{L}_{i, \boldsymbol{\pi}}$.

\begin{lemma}
\label{lInitialState}
We have
$$
\mathrm{Sel}_{\mathcal{L}_{0, \boldsymbol{\pi}}}(G_K, E[2]) = \mathrm{Sel}^2(E/K).
$$
\end{lemma}

\begin{proof}
Unwind the definitions and then use Lemma \ref{lEasy}$(b)$ for $v \not \in T$.
\end{proof}

For $v \in \Omega_K$ and $x \in \frac{K_v^\ast}{K_v^{\ast 2}}$, we have a natural isomorphism $E^x[2] \to E[2]$ as $G_{K_v}$-modules. Thus, after applying this natural identification, the quadratic forms $q_{E^x,v}$ and $q_{E,v}$ have the same cohomology group as source. The following establishes that these forms actually coincide. 

\begin{proposition}[{\cite[Lemma 5.2]{KMR}}]
\label{quadratic form is invariant under twisting}
Let $v \in \Omega_K$ and $x \in \frac{K_v^\ast}{K_v^{\ast 2}}$. Then $q_{E, v} = q_{E^x, v}$.
\end{proposition}

We will also use the following version of Poitou--Tate duality. We write $M^\vee = \Hom(M, \overline{K}^\ast)$ for the dual module, and we write $\mathcal{L}^\vee$ for the dual Selmer structure.

\begin{theorem}[{\cite[Theorem 2.3.4]{MR0}}]
\label{tPT}
Let $\mathcal{L}$ and $\mathcal{F}$ be Selmer structures such that $\mathcal{L}_v \subseteq \mathcal{F}_v$ for all $v \in \Omega_K$. Let $S$ be a finite set of places such that $\mathcal{L}_v = \mathcal{F}_v$ for all $v \not \in S$. Then the images of $\mathrm{Sel}_\mathcal{F}(G_K, M)$ and $\mathrm{Sel}_{\mathcal{L}^\vee}(G_K, M^\vee)$ are orthogonal complements under the pairing
\[
\bigoplus_{v \in S} \mathcal{F}_v/\mathcal{L}_v \times \bigoplus_{v \in S} \mathcal{L}_v^\vee/\mathcal{F}_v^\vee \rightarrow \Q/\Z
\]
given by the sum of the local Tate pairings over $S$.
\end{theorem}

\begin{proposition} 
\label{prop: relaxed Selmer image maximal isotropic}
Let $0 \leq i \leq n - 1$. Then $A_{v_{i + 1}}$ is a maximal isotropic subspace of the quadratic space $(H^1(G_{K_{v_{i + 1}}}, E[2]), q_{E,v_{i + 1}})$.
\end{proposition}

\begin{proof}
We apply Theorem \ref{tPT} with $\mathcal{F} := \mathcal{L}_{i, \boldsymbol{\pi}}'$ and $\mathcal{L} := (\mathcal{L}_v)_{v \in \Omega_K}$ given by
$$
\mathcal{L}_v :=
\begin{cases}
\mathcal{L}_{i, \boldsymbol{\pi}, v}' &\text{if } v \neq v_{i + 1}, \\
0 &\text{if } v = v_{i + 1}.
\end{cases}
$$
Note that $\mathcal{L}^\vee = \mathcal{F}$: this follows from Proposition \ref{prop: loc.glob.quad.form}$(a)$ for all $v$ different from $v_{i + 1}$ and holds by construction at $v = v_{i + 1}$. Then Theorem \ref{tPT} gives that $A_{v_{i + 1}}$ equals its own orthogonal complement under $b_{E, v_{i + 1}}$.

It now suffices to show that $q_{E, v_{i + 1}}$ vanishes on $A_{v_{i + 1}}$. We claim that
\begin{equation}
\label{eLocalqClaim}
q_{E, v}(\mathcal{L}_{i, \boldsymbol{\pi}, v}') = 0
\end{equation}
for all $v \in \Omega_K - \{v_{i + 1}\}$. If $v \in T$, this follows from the definition of $\mathcal{L}_{i, \boldsymbol{\pi}, v}'$ and Proposition \ref{prop: loc.glob.quad.form}$(a)$. If $v \not \in T \cup \{v_1, \dots, v_{i + 1}\}$, then the claim follows once more by the definition of $\mathcal{L}_{i, \boldsymbol{\pi}, v}'$, Lemma \ref{lEasy}$(b)$ and Proposition \ref{prop: loc.glob.quad.form}$(a)$. For $v \in \{v_1, \dots, v_i\}$, Lemma \ref{lEasy}$(c)$ gives
$$
\mathcal{L}_{i, \boldsymbol{\pi}, v}' = \delta \left(\frac{E^{\pi_i}(K_v)}{2E^{\pi_i}(K_v)}\right).
$$
Therefore the desired conclusion follows immediately from Proposition \ref{prop: loc.glob.quad.form}$(a)$ (applied to the quadratic twist $E^{\pi_i}$) and Proposition \ref{quadratic form is invariant under twisting}. We have now proven the claim \eqref{eLocalqClaim}.

Let $c \in \mathrm{Sel}_{\mathcal{L}_{i, \boldsymbol{\pi}}'}(G_K, E[2])$. Invoking Proposition \ref{prop: loc.glob.quad.form}$(b)$ and equation \eqref{eLocalqClaim}, we conclude that $q_{E, K}(c)$ is an element of $\text{Br}(K)[2]$ that restricts trivially at all $v$ different from $v_{i + 1}$. Therefore Hilbert reciprocity shows that $q_{E, K}(c)$ restricts trivially also at $v_{i + 1}$. Applying once more Proposition \ref{prop: loc.glob.quad.form}$(b)$, we conclude that $q_{E, v_{i + 1}}(\res_{v_{i + 1}}(c)) = 0$, as desired. 
\end{proof}

We deduce from Lemma \ref{lEasy}$(a)$ and Proposition \ref{prop: relaxed Selmer image maximal isotropic} that
\begin{align}
\label{eMaxIsotropic}
\dim_{\FF_2} A_{v_{i + 1}} = 2.
\end{align}
The additional quadratic structure was heavily used in the work of Klagsbrun--Mazur--Rubin \cite{KMR, KMR2}. The following key fact follows also from \cite[Proposition 7.2(iii)]{KMR2}, but we opt for a self-contained proof for the reader's convenience. 

\begin{lemma}
\label{lSelmerChange}
Let $0 \leq i \leq n - 1$. Then we have
$$
n_i =
\begin{cases}
2 &\textup{if } \dim_{\mathbb{F}_2} \textup{res}_{v_{i + 1}}(\mathrm{Sel}_{\mathcal{L}_{i, \boldsymbol{\pi}}}(G_K, E[2])) = 0 \textup{ and } A_{v_{i + 1}} = \mathcal{L}_{i + 1, \boldsymbol{\pi}, v_{i + 1}} \\
-2 &\textup{if } \dim_{\mathbb{F}_2} \textup{res}_{v_{i + 1}}(\mathrm{Sel}_{\mathcal{L}_{i, \boldsymbol{\pi}}}(G_K, E[2])) = 2 \\
0 &\textup{otherwise.}
\end{cases}
$$
\end{lemma}

\begin{proof}
Write $\overline{\res_{v_{i + 1}}}$ for the composition of $\res_{v_{i + 1}}$ with the natural quotient map going from $H^1(G_{K_{v_{i + 1}}}, E[2])$ to $H^1(G_{K_{v_{i + 1}}}, E[2])/\mathcal{L}_{i + 1, \boldsymbol{\pi}, v_{i + 1}}$. Let us start by recording the observation
\begin{equation}
\label{eSelmerTrick}
\mathrm{Sel}_{\mathcal{L}_{i + 1, \boldsymbol{\pi}}}(G_K, E[2]) = \ker\left(\mathrm{Sel}_{\mathcal{L}_{i, \boldsymbol{\pi}}'}(G_K, E[2]) \xrightarrow{\overline{\res_{v_{i + 1}}}} \frac{H^1(G_{K_{v_{i + 1}}}, E[2])}{\mathcal{L}_{i + 1, \boldsymbol{\pi}, v_{i + 1}}}\right)
\end{equation}
that we will use repeatedly.

Recalling that $A_{v_{i + 1}} = \res_{v_{i + 1}} \mathrm{Sel}_{\mathcal{L}_{i, \boldsymbol{\pi}}'}(G_K, E[2])$ by definition (see \eqref{eRelaxedA}), we also record the observations
\begin{alignat}{2}
&\res_{v_{i + 1}}(\mathrm{Sel}_{\mathcal{L}_{i, \boldsymbol{\pi}}}(G_K, E[2])) &&= A_{v_{i + 1}} \cap H^1_{\text{nr}}(G_{K_{v_{i + 1}}}, E[2]) \label{eResi} \\
&\res_{v_{i + 1}}(\mathrm{Sel}_{\mathcal{L}_{i + 1, \boldsymbol{\pi}}}(G_K, E[2])) &&= A_{v_{i + 1}} \cap \mathcal{L}_{i + 1, \boldsymbol{\pi}, v_{i + 1}} \label{eResi1},
\end{alignat}
which follow directly from unwinding the definitions of the various Selmer structures.

By first using equation \eqref{eResi} and then equation \eqref{eMaxIsotropic}, we obtain
$$
\dim_{\mathbb{F}_2} \res_{v_{i + 1}}(\mathrm{Sel}_{\mathcal{L}_{i, \boldsymbol{\pi}}}(G_K, E[2])) \leq \dim_{\mathbb{F}_2} A_{v_{i + 1}} = 2.
$$
We now examine the three possible cases for the value of this dimension. Suppose first that
\begin{equation}
\label{eDim2}
\dim_{\mathbb{F}_2} \res_{v_{i + 1}}(\mathrm{Sel}_{\mathcal{L}_{i, \boldsymbol{\pi}}}(G_K, E[2])) = 2.
\end{equation}
Since $\dim_{\FF_2} H^1_{\text{nr}}(G_{K_{v_{i + 1}}}, E[2]) = 2 = \dim_{\FF_2} A_{v_{i + 1}}$ by the second part of Lemma \ref{lEasy}$(a)$ and equation \eqref{eMaxIsotropic}, equation \eqref{eDim2} and equation \eqref{eResi} imply the equalities
\begin{equation}
\label{eAllEqual}
\res_{v_{i + 1}}(\mathrm{Sel}_{\mathcal{L}_{i, \boldsymbol{\pi}}}(G_K, E[2])) = A_{v_{i + 1}} = H^1_{\text{nr}}(G_{K_{v_{i + 1}}}, E[2]).
\end{equation}
In particular, the first equality in \eqref{eAllEqual} combined with the definition \eqref{eRelaxedA} implies the equality $\mathrm{Sel}_{\mathcal{L}_{i, \boldsymbol{\pi}}}(G_K, E[2]) = \mathrm{Sel}_{\mathcal{L}_{i, \boldsymbol{\pi}}'}(G_K, E[2])$. Therefore equation \eqref{eSelmerTrick} yields
\begin{equation}
\label{eSubspacedim2}
\mathrm{Sel}_{\mathcal{L}_{i + 1, \boldsymbol{\pi}}}(G_K, E[2]) = \ker\left(\mathrm{Sel}_{\mathcal{L}_{i, \boldsymbol{\pi}}}(G_K, E[2]) \xrightarrow{\overline{\res_{v_{i + 1}}}} \frac{H^1(G_{K_{v_{i + 1}}}, E[2])}{\mathcal{L}_{i + 1, \boldsymbol{\pi}, v_{i + 1}}}\right).
\end{equation}
Moreover, by equation \eqref{eAllEqual}, the image of $\mathrm{Sel}_{\mathcal{L}_{i, \boldsymbol{\pi}}}(G_K, E[2])$ under $\overline{\res_{v_{i + 1}}}$ is
$$
\frac{\mathcal{L}_{i + 1, \boldsymbol{\pi}, v_{i + 1}} A_{v_{i + 1}}}{\mathcal{L}_{i + 1, \boldsymbol{\pi}, v_{i + 1}}}.
$$
Since $\mathcal{L}_{i +1, \boldsymbol{\pi}, v_{i + 1}}$ intersects $H^1_{\text{nr}}(G_{K_{v_{i + 1}}}, E[2])$ trivially by Lemma \ref{lEasy}$(d)$ and since we have $A_{v_{i + 1}} = H^1_{\text{nr}}(G_{K_{v_{i + 1}}}, E[2])$ by equation \eqref{eAllEqual}, it follows from equation \eqref{eMaxIsotropic} that $n_i = -2$.

Suppose now that 
\begin{equation}
\label{eDim0}
\dim_{\mathbb{F}_2} \res_{v_{i + 1}}(\mathrm{Sel}_{\mathcal{L}_{i, \boldsymbol{\pi}}}(G_K, E[2])) = 0.
\end{equation}
We will distinguish two further cases, namely $A_{v_{i + 1}} = \mathcal{L}_{i + 1, \boldsymbol{\pi}, v_{i + 1}}$ or $A_{v_{i + 1}} \neq \mathcal{L}_{i + 1, \boldsymbol{\pi}, v_{i + 1}}$. 

Suppose that $A_{v_{i + 1}} = \mathcal{L}_{i + 1, \boldsymbol{\pi}, v_{i + 1}}$. Then equation \eqref{eSelmerTrick} shows that
\begin{equation}
\label{eSelJump}
\mathrm{Sel}_{\mathcal{L}_{i + 1, \boldsymbol{\pi}}}(G_K, E[2]) = \mathrm{Sel}_{\mathcal{L}_{i, \boldsymbol{\pi}}'}(G_K, E[2]).
\end{equation}
Since $A_{v_{i + 1}} = \mathcal{L}_{i + 1, \boldsymbol{\pi}, v_{i + 1}}$ intersects $H^1_{\text{nr}}(G_{K_{v_{i + 1}}}, E[2])$ trivially by Lemma \ref{lEasy}$(d)$, we deduce from equation \eqref{eResi} that $\mathrm{Sel}_{\mathcal{L}_{i, \boldsymbol{\pi}}}(G_K, E[2])$ equals the subspace of $\mathrm{Sel}_{\mathcal{L}_{i, \boldsymbol{\pi}}'}(G_K, E[2])$ that restricts trivially at $v_{i + 1}$. Hence we have $n_i = 2$ by equations \eqref{eMaxIsotropic} and \eqref{eSelJump}.

Now suppose that $A_{v_{i + 1}} \neq \mathcal{L}_{i + 1, \boldsymbol{\pi}, v_{i + 1}}$. We have
\begin{equation}
\label{eAvunr}
A_{v_{i + 1}} \cap H_{\text{nr}}^1(G_{K_{v_{i + 1}}}, E[2]) = \{0\}
\end{equation}
by equations \eqref{eResi} and \eqref{eDim0}. Observe that $\mathcal{L}_{i + 1, \boldsymbol{\pi}, v_{i + 1}}$ intersects $H^1_{\text{nr}}(G_{K_{v_{i + 1}}}, E[2])$ trivially by Lemma \ref{lEasy}$(d)$, and observe moreover that $A_{v_{i + 1}}$, $\mathcal{L}_{i + 1, \boldsymbol{\pi}, v_{i + 1}}$ and $H^1_{\text{nr}}(G_{K_{v_{i + 1}}}, E[2])$ are all maximal isotropic subspaces of $H^1(G_{K_{v_{i + 1}}}, E[2])$ (for $A_{v_{i + 1}}$, see Proposition \ref{prop: relaxed Selmer image maximal isotropic}, while the other two spaces are of the shape $\delta(E^x(K_v)/2E^x(K_v))$ for some quadratic twist $E^x$ of $E$ by Lemma \ref{lEasy} and are thus maximal isotropic by Proposition \ref{prop: loc.glob.quad.form}$(a)$ and Proposition \ref{quadratic form is invariant under twisting}). Combining this with equation \eqref{eAvunr}, we see that \cite[Proposition 2.4]{KMR} yields
$$
\dim_{\mathbb{F}_2} A_{v_{i + 1}} \cap \mathcal{L}_{i + 1, \boldsymbol{\pi}, v_{i + 1}} \equiv 0 \bmod 2.
$$
However, the left hand side is strictly smaller than $2$ since $A_{v_{i + 1}} \neq \mathcal{L}_{i + 1, \boldsymbol{\pi}, v_{i + 1}}$ by assumption (and both spaces are $2$-dimensional). Therefore the left hand side has to be equal to $0$, i.e. 
\begin{equation}
\label{eAvL}
A_{v_{i + 1}} \cap \mathcal{L}_{i + 1, \boldsymbol{\pi}, v_{i + 1}} = \{0\}.
\end{equation}
We conclude that
$$
\mathrm{Sel}_{\mathcal{L}_{i, \boldsymbol{\pi}}}(G_K, E[2]) = \ker\left(\mathrm{Sel}_{\mathcal{L}_{i, \boldsymbol{\pi}}'}(G_K, E[2]) \xrightarrow{\res_{v_{i + 1}}} H^1(G_{K_{v_{i + 1}}}, E[2])\right) = \mathrm{Sel}_{\mathcal{L}_{i + 1, \boldsymbol{\pi}}}(G_K, E[2]),
$$
where the first equality follows from equations \eqref{eAvunr} and \eqref{eResi} and the second equality follows from equations \eqref{eAvL} and \eqref{eResi1}. In particular, we have $n_i = 0$ as desired.

It remains to consider the case
$$
\dim_{\mathbb{F}_2} \res_{v_{i + 1}}(\mathrm{Sel}_{\mathcal{L}_{i, \boldsymbol{\pi}}}(G_K, E[2])) = 1.
$$
By equation \eqref{eResi}, this may be rephrased as
\begin{equation}
\label{eAvunr2}
\dim_{\mathbb{F}_2} A_{v_{i + 1}} \cap H^1_{\text{nr}}(G_{K_{v_{i + 1}}}, E[2]) = 1.
\end{equation}
Recall that $\mathcal{L}_{i + 1, \boldsymbol{\pi}, v_{i + 1}}$ intersects $H^1_{\text{nr}}(G_{K_{v_{i + 1}}}, E[2])$ trivially by Lemma \ref{lEasy}$(d)$. Combining this with equation \eqref{eAvunr2}, we obtain from \cite[Proposition 2.4]{KMR} the congruence
$$
\dim_{\mathbb{F}_2} A_{v_{i + 1}} \cap \mathcal{L}_{i + 1, \boldsymbol{\pi}, v_{i + 1}} \equiv 1 \bmod 2.
$$
Therefore the above dimension must be exactly equal to $1$. Thus, $\mathrm{Sel}_{\mathcal{L}_{i + 1, \boldsymbol{\pi}}}(G_K, E[2])$ equals the inverse image of the $1$-dimensional space $A_{v_{i + 1}} \cap \mathcal{L}_{i + 1, \boldsymbol{\pi}, v_{i + 1}}$ under the map
$$
\res_{v_{i + 1}}: \mathrm{Sel}_{\mathcal{L}_{i,\boldsymbol{\pi}}'}(G_K, E[2]) \twoheadrightarrow A_{v_{i + 1}}.
$$
At the same time $\mathrm{Sel}_{\mathcal{L}_{i, \boldsymbol{\pi}}}(G_K, E[2])$ equals the inverse image of $A_{v_{i + 1}} \cap H_{\text{nr}}^1(G_{K_{v_{i + 1}}}, E[2])$ under the same map, which is also $1$-dimensional by equation \eqref{eAvunr2}. Hence the two spaces $\mathrm{Sel}_{\mathcal{L}_{i,\boldsymbol{\pi}}}(G_K, E[2])$ and $\mathrm{Sel}_{\mathcal{L}_{i + 1,\boldsymbol{\pi}}}(G_K, E[2])$ are both subspaces in $\mathrm{Sel}_{\mathcal{L}_{i,\boldsymbol{\pi}}'}(G_K, E[2])$ of codimension $1$ and hence they must have the same dimension. We conclude that $n_i = 0$, as desired.
\end{proof}

\begin{remark}
\label{rSelmerChangeUpgrade}
In some of our applications, it will be important to strengthen the conclusion of Lemma \ref{lSelmerChange} as follows:
\begin{enumerate}
\item[(a)] the condition $A_{v_{i + 1}} = \mathcal{L}_{i + 1, \boldsymbol{\pi}, v_{i + 1}}$ implies $\dim_{\mathbb{F}_2} \res_{v_{i + 1}}(\mathrm{Sel}_{\mathcal{L}_{i, \boldsymbol{\pi}}}(G_K, E[2])) = 0$. Indeed, this follows from the fact that $\res_{v_{i + 1}}(\mathrm{Sel}_{\mathcal{L}_{i, \boldsymbol{\pi}}}(G_K, E[2]))$ is exactly 
$$
A_{v_{i + 1}} \cap H^1_{\textup{nr}}(G_{K_{v_{i + 1}}}, E[2]) = \mathcal{L}_{i + 1, \boldsymbol{\pi}, v_{i + 1}} \cap H^1_{\textup{nr}}(G_{K_{v_{i + 1}}}, E[2]) = \{0\}
$$
by Lemma \ref{lEasy}$(d)$;
\item[(b)] $n_i = -2$ implies that $\mathrm{Sel}_{\mathcal{L}_{i + 1, \boldsymbol{\pi}}}(G_K, E[2]) \subseteq \mathrm{Sel}_{\mathcal{L}_{i, \boldsymbol{\pi}}}(G_K, E[2])$ is a subspace of codimension $2$ (indeed, combine $n_i = -2$ with equation \eqref{eSubspacedim2});
\item[(c)] since $\textup{res}_{v_{i + 1}}(\mathrm{Sel}_{\mathcal{L}_{i, \boldsymbol{\pi}}}(G_K, E[2])) \subseteq H^1_{\textup{nr}}(G_{K_{v_{i + 1}}}, E[2])$ by definition and since moreover $\dim_{\FF_2} H^1_{\textup{nr}}(G_{K_{v_{i + 1}}}, E[2]) = 2$ by the second part of Lemma \ref{lEasy}$(a)$, the condition 
$$
\dim_{\mathbb{F}_2} \textup{res}_{v_{i + 1}}(\mathrm{Sel}_{\mathcal{L}_{i, \boldsymbol{\pi}}}(G_K, E[2])) = 2
$$
is equivalent to the existence of two linearly independent elements inside the restriction $\textup{res}_{v_{i + 1}}(\mathrm{Sel}_{\mathcal{L}_{i, \boldsymbol{\pi}}}(G_K, E[2]))$.
\end{enumerate}
\end{remark}

The final step in the transition process is completely predictable because of the $2$-torsion elements of $E$.

\begin{lemma}
\label{lFinalStep}
Let $d \in K^\ast/K^{\ast 2}$. Assume that
\begin{equation}
\label{eTwistLocalTrivial}
d \in K_v^{\ast 2} \quad \textup{ for all } v \in T.
\end{equation}
Let $v_1, \dots, v_r$ be the set of places $v \in \Omega_K$ such that $v(d)$ is odd; note that $v_1, \dots, v_r \not \in T$ by \eqref{eTwistLocalTrivial}. We assume that $r \geq 1$. Then we have
$$
\dim_{\mathbb{F}_2} \mathrm{Sel}^2(E^d/K) = 2 + \dim_{\mathbb{F}_2} \mathrm{Sel}_{\mathcal{L}_{r - 1, (d, \dots, d)}}(G_K, E[2]).
$$
\end{lemma}

\begin{proof}
We define the Selmer structure $\mathcal{L}_d = (\mathcal{L}_{d, v})_{v \in \Omega_K}$ by
$$
\mathcal{L}_{d, v} := \delta(E^d(K_v)).
$$
Recall that $\mathrm{Sel}^2(E^d/K) = \mathrm{Sel}_{\mathcal{L}_d}(G_K, E[2])$ by definition. We remark that 
\begin{align}
\label{eEd}
\mathcal{L}_{d, v} = \delta(E^d(K_v)) = 
\begin{cases}
\delta(E(K_v)) &\text{if } v \in T \\
\langle (\alpha \beta, d \alpha), (-d \alpha, -\alpha \gamma) \rangle &\text{if } v \in \{v_1, \dots, v_r\} \\
H^1_{\text{nr}}(G_{K_v}, E[2]) &\text{if } v \not \in T \cup \{v_1, \dots, v_r\}.
\end{cases}
\end{align}
The first row follows from the isomorphism $E \cong E^d$ locally at $K_v$, which is a consequence of our assumption \eqref{eTwistLocalTrivial}. The second and third row follow from respectively Lemma \ref{lEasy}$(c)$ and Lemma \ref{lEasy}$(b)$.

Set $\boldsymbol{\pi} := (d, \dots, d)$. By comparing equation \eqref{eSelmerStructure} with equation \eqref{eEd}, it follows that $\mathcal{L}_{d, v} = \mathcal{L}_{r, \boldsymbol{\pi}, v}$ for all $v \in \Omega_K$. Therefore we certainly have
$$
\mathrm{Sel}^2(E^d/K) = \mathrm{Sel}_{\mathcal{L}_{r, \boldsymbol{\pi}}}(G_K, E[2]),
$$
so it suffices to show that $n_{r - 1} = 2$. 

We know that $A_{v_r}$ has dimension $2$ by equation \eqref{eMaxIsotropic}. However, it also contains (the restrictions of) $(\alpha \beta, d \alpha)$ and $(-d \alpha, -\alpha \gamma)$ as these are in $\mathrm{Sel}^2(E^d/K)$ by Lemma \ref{lemma: delta map formula} and hence in $\mathrm{Sel}_{\mathcal{L}_{r, \boldsymbol{\pi}}}(G_K, E[2])$. Since $r \geq 1$, it follows that $\res_{v_r} (\alpha \beta, d \alpha)$ and $\res_{v_r} (-d \alpha, -\alpha \gamma)$ are linearly independent by inspecting valuations at $v_r$ (recall that $v_r \not \in T$ by \eqref{eTwistLocalTrivial}). Therefore $A_{v_r}$ is equal to the span of $(\alpha \beta, d \alpha)$ and $(-d \alpha, -\alpha \gamma)$, which is $\mathcal{L}_{r, \boldsymbol{\pi}, v_r}$. Thus we obtain the equality $n_{r - 1} = 2$ from Remark \ref{rSelmerChangeUpgrade}$(a)$ and Lemma \ref{lSelmerChange}, as desired.
\end{proof}
\section{Proof of Theorem \ref{tRelative}}
\label{sGT}
In this section we combine the material from Section \ref{sBackground} and \ref{sGen} to prove Theorem \ref{tRelative}. For the remainder of the paper, we fix a number field $K$ satisfying the conditions of Theorem \ref{tRelative}, so $K$ has at least $32$ real places and we set $L := K(i)$. 

\subsection{First steps}
Let $a_1$, $a_2$, $a_3$, $E$, $\sigma_1, \dots, \sigma_{24}$ be as in Lemma \ref{lStartingCurve}. Let $\alpha$, $\beta$, $\gamma$ be as in equation \eqref{eabc}. We will twist $E$. Define 
$$
T := \{v \in \Omega_K : v \mid 6 \Delta_E \infty \text{ or } v(\alpha) \neq 0 \text{ or } v(\beta) \neq 0 \text{ or } v(\gamma) \neq 0\}.
$$
In the next two subsections, we will find $t \in K^\ast/K^{\ast 2}$ satisfying
\begin{align}
\label{eKey}
\dim_{\mathbb{F}_2} \, \mathrm{Sel}^2(E^t/K) = 2, \quad \mathrm{rk} \, E^{-t}(K) > 0.
\end{align}
Because $E^t$ has full rational $2$-torsion, we have the inequality 
$$
\mathrm{rk} \, E^t(K) + 2 \leq \dim_{\mathbb{F}_2} \, \mathrm{Sel}^2(E^t/K).
$$
Hence it follows from the first part of equation \eqref{eKey} that $\mathrm{rk} \, E^t(K) = 0$. Combining this with the second part of equation \eqref{eKey} and Lemma \ref{lRankDecomposition} (applied to $E^{-t}$ and $d = -1$) we deduce that $\mathrm{rk} \, E^{-t}(K) = \mathrm{rk} \, E^{-t}(L) > 0$. Hence Theorem \ref{tRelative} holds with $E^{-t}$ being our witness.

We will establish \eqref{eKey} with two reduction steps. The first reduction step gives sufficient conditions on $t$ (see Definition \ref{dSui}) to ensure that 
$$
\dim_{\mathbb{F}_2} \, \mathrm{Sel}^2(E^t/K) = 2. 
$$
The second reduction step finds sufficient conditions on $\kappa$ (see Definition \ref{dAux} for a precise statement) that will allow us to construct $t$ through additive combinatorics. Specifically, $t$ will take the shape
$$
t = \kappa q_1 q_2 q_3 q_3 q_4,
$$
where $\kappa$ is an auxiliary twist to arrange the Selmer group favorably, and $q_1, q_2, q_3, q_4$ are constructed through additive combinatorics to ensure that equation \eqref{eKey} holds. Finally, we finish the proof by showing that an auxiliary twist $\kappa$ exists.

% Let $t \in K^{\ast}/K^{\ast 2}$ and let $v_1, \dots, v_r$ be the places ramified in $\psi_t$ outside of $T$. For each integer $0 \leq i \leq r$, we define the Selmer structure $\mathcal{L}_{i, t}$ to be $\mathcal{L}_{i, \boldsymbol{\pi}}$ from Section \ref{sGen}, where $\mathrm{res}_{v_j}(\psi_t) = \psi_{\pi_j}$. We will identify $\frac{1}{2}\Z/\Z \cong \mathbb{F}_2$, so that our invariant maps will naturally land in $\mathbb{F}_2$. 

\begin{mydef}
\label{dSui}
We say that an element $t \in K^\ast/K^{\ast 2}$ is a suitable twist if there exists $\kappa \in K^\ast/K^{\ast 2}$ and prime elements $q_1, q_2, q_3, q_4$ such that, denoting by $\mathfrak{p}_1, \dots, \mathfrak{p}_s$ the primes where $\kappa$ has odd valuation and putting $T' := T \cup \{\mathfrak{p}_1, \dots, \mathfrak{p}_s\}$, then
\begin{enumerate}
\item[$(P1)$] $t \in K_v^{\ast 2}$ for all places $v \in T$;
\item[$(P2)$] we have $t = \kappa q_1 q_2 q_3 q_4$ and the elements $q_1, q_2, q_3, q_4$ are pairwise coprime and also coprime to $T'$;
\item[($P3)$] there exists a basis
$$
(z_1, z_2), (z_3, z_4), (z_5, z_6), (z_7, z_8), (z_9, z_{10}), (z_{11}, z_{12})
$$
of $\mathrm{Sel}_{\mathcal{L}_{s, (t, \dots, t)}}(G_K, E[2])$ such that
\begin{align*}
&\prod_{v \in T'} (z_1, q_1)_v = \prod_{v \in T'} (z_4, q_1)_v = -1 \\
&\prod_{v \in T'} (z_5, q_2)_v = \prod_{v \in T'} (z_8, q_2)_v = -1 \\
&\prod_{v \in T'} (z_9, q_3)_v = \prod_{v \in T'} (z_{12}, q_3)_v = -1 \\
&\prod_{v \in T'} (z_i, q_1)_v = 1 \quad \textup{ for } i \in \{1, \dots, 12\} - \{1, 4\} \\
&\prod_{v \in T'} (z_i, q_2)_v = 1 \quad \textup{ for } i \in \{5, \dots, 12\} - \{5, 8\} \\
&\prod_{v \in T'} (z_i, q_3)_v = 1 \quad \textup{ for } i \in \{9, \dots, 12\} - \{9, 12\};
\end{align*}
\item[$(P4)$] we have $\mathrm{rk} \, E^{-t}(K) > 0$.
\end{enumerate}
\end{mydef}

Property $(P3)$ serves two distinct roles. Firstly, we will use this property to understand how $q_1, q_2, q_3$ behave in the transition process from Lemma \ref{lSelmerChange}. In fact, we will show that the dimension drops by $2$ in each case. Secondly, $(P3)$ imposes merely congruences on $q_1, q_2, q_3$ modulo $T'$. This is very convenient, as Theorem \ref{tKai2} does allow us to prescribe such congruence conditions, but does not allow for control on mutual Legendre symbols $(q_i/q_j)$. Condition $(P3)$ carefully sets up a situation where we can understand the change of Selmer rank without having knowledge of these Legendre symbols.

\begin{theorem}
\label{tFirstReduction}
Let $t \in K^{\ast}/K^{\ast 2}$ be suitable. Then $t$ satisfies equation \eqref{eKey}.
\end{theorem}

\begin{proof}
Because $t$ is suitable, $(P4)$ holds. This is exactly the second part of equation \eqref{eKey}, and thus we need to show that
\begin{equation}
\label{eKey2}
\dim_{\mathbb{F}_2} \, \mathrm{Sel}^2(E^t/K) = 2.
\end{equation}
By $(P1)$ and $(P2)$, we have that 
$$
\{v \in \Omega_K : v \text{ finite}, v(t) \equiv 1 \bmod 2\} = \{\mathfrak{p}_1, \dots, \mathfrak{p}_s, (q_1), (q_2), (q_3), (q_4)\}.
$$
We claim that
\begin{align}
\label{eIntermediate}
\mathrm{Sel}_{\mathcal{L}_{s + 1, (t, \dots, t)}}(G_K, E[2]) = \langle (z_5, z_6), (z_7, z_8), (z_9, z_{10}), (z_{11}, z_{12}) \rangle.
\end{align}
To this end, we apply Hilbert reciprocity to obtain the relations
$$
\prod_{v \in \Omega_K} (z_i, q_1)_v = 1
$$
for all $i \in \{1, \dots, 12\}$. Using that $v(z_i)$ is even for all $v \not \in T'$, we deduce that $(z_i, q_1)_v = 1$ for all $v \not \in T' \cup \{(q_1)\}$. Therefore $(P3)$ reveals the identity
\begin{equation}
\label{eReductionzi}
(z_i, q_1)_{(q_1)} = 
\begin{cases}
-1 &\text{if } i \in \{1, 4\} \\
1 &\text{otherwise.}
\end{cases}
\end{equation}
Now note that $(z_i, q_1)_{(q_1)}$ detects precisely whether $\mathrm{res}_{(q_1)}(z_i)$ is the trivial class or the unique non-trivial unramified class $\epsilon$. Combining this with equation \eqref{eReductionzi} we conclude that
\begin{equation}
\label{eReductionWrittenOut}
\mathrm{res}_{(q_1)}((z_1, z_2)) = (\epsilon, 1), \quad \mathrm{res}_{(q_1)}((z_3, z_4)) = (1, \epsilon), \quad \mathrm{res}_{(q_1)}((z_{2i + 1}, z_{2i + 2})) = (1, 1)
\end{equation}
for all $i \in \{2, 3, 4, 5\}$. Hence Lemma \ref{lSelmerChange} implies that $n_s = -2$. Moreover, from equation \eqref{eReductionWrittenOut}, we visibly have that $(z_5, z_6)$, $(z_7, z_8)$, $(z_9, z_{10})$ and $(z_{11}, z_{12})$ are in $\mathrm{Sel}_{\mathcal{L}_{s + 1, (t, \dots, t)}}(G_K, E[2])$. Since these elements are linearly independent, we deduce equation \eqref{eIntermediate}.

Proceeding in this way, one successively shows that
$$
\mathrm{Sel}_{\mathcal{L}_{s + 2, (t, \dots, t)}}(G_K, E[2]) = \langle (z_9, z_{10}), (z_{11}, z_{12}) \rangle
$$
and that
\begin{equation}
\label{eEndStep}
\dim_{\mathbb{F}_2} \, \mathrm{Sel}_{\mathcal{L}_{s + 3, (t, \dots, t)}}(G_K, E[2]) = 0
\end{equation}
for all suitable $t$. 

Note that equation \eqref{eTwistLocalTrivial} holds thanks to our assumption $(P1)$. Therefore the hypotheses of Lemma \ref{lFinalStep} are satisfied (note that the condition $r \geq 1$ is clearly satisfied in our case). Then equation \eqref{eKey2} is an immediate consequence of equation \eqref{eEndStep} and Lemma \ref{lFinalStep}.
\end{proof}

Thanks to Theorem \ref{tFirstReduction}, it remains to establish the existence of a suitable $t$.

\subsection{Application of additive combinatorics}
We will now construct a suitable $t$ in two steps. First, we shall assume that a certain auxiliary twist $\kappa$ has certain good properties. Under this assumption, we will apply Kai's work \cite{Kai} to find such $t$. The construction of $\kappa$ is delicate and relies on the material in Section \ref{sGen} in an essential way. This construction will occupy the entirety of Subsection \ref{ssFinal}.

By the second part of Lemma \ref{lStartingCurve}, we fix six infinite places $\tau_1, \dots, \tau_6$ of $K$ among $\sigma_1, \dots, \sigma_{24}$ with
\begin{alignat}{2}
&\tau_1(a_2) < \tau_1(a_1) < \tau_1(a_3), \quad \quad % So tau_1(alpha beta) < 0, tau_1(alpha) > 0
&&\tau_2(a_1) < \tau_2(a_2) < \tau_2(a_3) \label{eTau12} \\ % So \tau_2(\alpha \beta) > 0, \tau_2(\alpha) < 0
&\tau_3(a_2) < \tau_3(a_3) < \tau_3(a_1), \quad \quad % So \tau_3(\alpha \beta) > 0, \tau_3(\alpha) < 0
&&\tau_4(a_1) < \tau_4(a_3) < \tau_4(a_2) \label{eTau34} \\ % So \tau_4(\alpha) > 0, \tau_4(\alpha \gamma) > 0
&\tau_5(a_2) < \tau_5(a_1) < \tau_5(a_3), \quad \quad % So tau_5(alpha beta) > 0, tau_5(alpha) < 0
&&\tau_6(a_1) < \tau_6(a_2) < \tau_6(a_3). \label{eTau56} % So \tau_6(\alpha \beta) > 0, \tau_6(\alpha) < 0
\end{alignat}

\begin{mydef}
\label{dAux}
We say that $\kappa \in K^{\ast}/K^{\ast 2}$ is an auxiliary twist if
\begin{enumerate}
\item[$(K1)$] $\kappa$ has even valuation at all finite places in $T$, and is non-trivial locally at $\tau_4, \tau_5, \tau_6$, but trivial at $\tau_1, \tau_2, \tau_3$;
\item[$(K2)$] writing $\mathfrak{p}_1, \dots, \mathfrak{p}_s$ for the primes where $\kappa$ has odd valuation, there exists $\boldsymbol{\pi} = (\pi_i)_{1 \leq i \leq s}$ such that $\mathrm{Sel}_{\mathcal{L}_{s, \boldsymbol{\pi}}}(G_K, E[2])$ has a basis 
$$
(z_1, z_2), (z_3, z_4), (z_5, z_6), (z_7, z_8), (z_9, z_{10}), (z_{11}, z_{12})
$$
such that the signs of $z_i$ at $\tau_1, \dots, \tau_6$ are given by the table
\begin{equation}
\label{eInfMat1}
\begin{array}{| c | c | c | c | c | c | c |}
\hline & \tau_1 & \tau_2 & \tau_3 & \tau_4 & \tau_5 & \tau_6 \\
\hline \text{im}(\delta) & (-, +) & (+, -) & (-, +) & (+, -) & (-, +) & (+, -) \\
\hline z_1    & - & + & + & + & + & + \\
\hline z_2    & + & + & + & + & + & + \\
\hline z_3    & + & + & + & + & + & + \\
\hline z_4    & + & - & + & + & + & + \\
\hline z_5    & + & + & - & + & + & + \\
\hline z_6    & + & + & + & + & + & + \\
\hline z_7    & + & + & + & + & + & + \\
\hline z_8    & + & + & + & - & + & + \\
\hline z_9    & + & + & + & + & - & + \\
\hline z_{10} & + & + & + & + & + & + \\
\hline z_{11} & + & + & + & + & + & + \\
\hline z_{12} & + & + & + & + & + & - \\ \hline
\end{array}
\end{equation}
\end{enumerate}
\end{mydef}

\noindent For the reader's convenience, we have added the image $\delta(E(K_{\tau_j}))$ of the local Kummer map at each real place $\tau_j$ in the second row of \eqref{eInfMat1}. This image can be directly computed by using the explicit formula from Lemma \ref{lemma: delta map formula}. In particular, one may check that the pair $(z_{2i + 1}, z_{2i + 2})$ satisfies the Selmer conditions at each $\tau_j$. 

We will now give some motivation for $(K2)$. Write $N$ for the product of all odd prime ideals in $T$. In order to construct suitable $t$, we shall soon choose the primes $q_1, q_2, q_3$ in such a way that 
$$
q_1 \equiv q_2 \equiv q_3 \equiv 1 \bmod 8 N \mathfrak{p}_1 \cdots \mathfrak{p}_s
$$
and that $q_1, q_2, q_3$ are positive at all real places except $\{\tau_1, \dots, \tau_6\}$. Then the sum over $v \in T'$ in $(P3)$ will simplify to a sum over just $v \in \{\tau_1, \dots, \tau_6\}$. With the further sign choices we shall make for $q_1, q_2, q_3$ at $\tau_1, \dots, \tau_6$ in \eqref{eInfMat2}, we will see that the table \eqref{eInfMat1} is constructed exactly to give back $(P3)$.

We are now ready to phrase our last reduction step, after which it remains to construct an auxiliary twist $\kappa$ in Subsection \ref{ssFinal}.

\begin{theorem}
\label{tSuitable}
Assume that there exists an auxiliary twist $\kappa$. Then there exists a suitable $t \in K^\ast/K^{\ast 2}$.
\end{theorem}

\begin{proof}
We will proceed in four steps, aiming to eventually apply Theorem \ref{tKai2}. In our first step, we construct four affine linear polynomials $L_1, L_2, L_3, L_4$. In the second step, we construct a convex region $\Omega \subseteq \mathbb{R}^{2n}$. In the third step we estimate the volume of $\Omega \cap [-H, H]^{2n}$. Finally, we apply Theorem \ref{tKai2} and verify the properties $(P1)$, $(P2)$, $(P3)$ and $(P4)$ in the fourth step.

\subsubsection*{Construction of affine linear polynomials}
Fix an auxiliary twist $\kappa$. Then, by $(K1)$ and strong approximation, there exists an integral representative of $\kappa$ coprime to $8N$. We also fix such an integral representative, which by abuse of notation we also denote by $\kappa$. 

Writing $h$ for the class number of $K$, we let $m$ be a generator of the ideal $8 N^h$. Write $\lambda \in O_K$ for an element coprime to $\kappa$ and the finite places in $T$ such that $\lambda$ reduces to the inverse of $\kappa$ modulo $8N$ and satisfies $\mathrm{res}_{\mathfrak{p}_i}(\lambda \kappa) = \pi_i$ for all $1 \leq i \leq s$. We define affine linear polynomials $L_1, L_2, L_3, L_4 \in O_K[X, Y]$ by
\begin{align*}
L_1(X, Y) &:= m^2 \kappa X + a_1 m^2 \kappa (m^2 \kappa Y + \lambda) + 1 \\
L_2(X, Y) &:= m^2 \kappa X + a_2 m^2 \kappa (m^2 \kappa Y + \lambda) + 1 \\ 
L_3(X, Y) &:= m^2 \kappa X + a_3 m^2 \kappa (m^2 \kappa Y + \lambda) + 1 \\
L_4(X, Y) &:= m^2 \kappa Y + \lambda.
\end{align*}
Since the $a_i$ are distinct, it follows immediately that the $L_i(X, Y)$ are pairwise $K$-linearly independent.

We identify $O_K$ with the group $\Z^n$ to reinterpret the $L_i$ as linear polynomials $\Z^{2n} \rightarrow O_K$. In order to show that the main term in Theorem \ref{tKai2} dominates, we need to check that the $p$-adic densities $\beta_p^{O_K}(L_1, L_2, L_3, L_4)$ there do not vanish. We will do this now.

\begin{lemma}
\label{lAdmissible}
For every prime $\mathfrak{p}$ of $O_K$, there exist $u, v \in O_K$ such that 
$$
L_1(u, v) L_2(u, v) L_3(u, v) L_4(u, v) \not \equiv 0 \bmod \mathfrak{p}.
$$
\end{lemma}

\begin{proof}
For primes $\mathfrak{p}$ dividing $m \kappa$, we may pick $u = v = 0$. It remains to verify the lemma for primes $\mathfrak{p} \nmid m \kappa$. In particular, $\mathfrak{p} \nmid 6$, so the residue field $O_K/\mathfrak{p}$ is $\mathbb{F}_q$ for some $q \geq 5$. Since $m^2 \kappa$ is a unit in $\mathbb{F}_q$, each $L_i$ cuts out a line in $\mathbb{F}_q \times \mathbb{F}_q$, leaving at least $q^2 - 4q > 0$ possibilities for $(u, v) \bmod \mathfrak{p}$.
\end{proof}

\subsubsection*{Construction of a convex region}
We will now construct a convex set $\Omega$ to which we apply Theorem \ref{tKai2}. We define $\Omega \subseteq \mathbb{R}^{2n}$ through the simultaneous conditions:
\begin{enumerate}
\item[$(I)$] for all real embeddings $\tau$ of $K$ not equal to $\tau_1, \dots, \tau_6$ and all $k \in \{1, 2, 3\}$
\begin{align}
\label{eEasyInf1}
\tau\left(L_k(X, Y)\right) > 0.
\end{align}
Here we note that $\tau \circ L_k$ is an affine $\Z$-linear map $\Z^{2n} \rightarrow \mathbb{R}$, which extends uniquely to an affine $\R$-linear map $\mathbb{R}^{2n} \rightarrow \mathbb{R}$. From now on we shall frequently extend linear maps $\Z^{2n} \rightarrow \mathbb{R}$ uniquely to $\mathbb{R}^{2n} \rightarrow \mathbb{R}$ without further mention;
\item[$(II)$] for all real embeddings $\tau$ of $K$ not equal to $\tau_1, \dots, \tau_6$
\begin{align}
\label{eEasyInf2}
\tau\left(L_4(X, Y)\right) \tau(\kappa) > 0;
\end{align}
\item[$(III)$] for $j \in \{1, \dots, 6\}$ and $k \in \{1, 2, 3\}$, the sign of $\tau_j\left(L_k(X, Y)\right)$ is precisely the entry $(k, j)$ in the table
\begin{equation}
\label{eInfMat2}
\begin{array}{| c | c | c | c | c | c | c |}
\hline & \tau_1 & \tau_2 & \tau_3 & \tau_4 & \tau_5 & \tau_6 \\
\hline L_1(X, Y)    & - & - & + & + & + & + \\
\hline L_2(X, Y)    & - & - & - & - & + & + \\
\hline L_3(X, Y)    & + & + & - & + & - & - \\ \hline
\end{array}
\end{equation}
\item[$(IV)$] for all $j \in \{1, \dots, 6\}$
$$
\tau_j\left(L_4(X, Y)\right) > 0.
$$
\end{enumerate}

\subsubsection*{Volume computation}
In order to show that the main term in Theorem \ref{tKai2} dominates, we also need to asymptotically estimate the volume of $\Omega \cap [-H, H]^{2n}$. We start with a general lemma on volumes. Write $\mu$ for the Lebesgue measure.

\begin{lemma}
\label{lRealVolume}
Let $\ell \geq k \geq 1$ be integers, let $A: \R^\ell \twoheadrightarrow \R^k$ be a surjective linear operator, and let $(c_1, \ldots, c_k) \in \R^k$. Then there exists a positive real number $C$ such that
$$
\mu\left(A^{-1}\left(\prod_{i = 1}^k (c_i, \infty)\right) \cap [-H, H]^\ell \right) \sim C H^\ell
$$
as $H \rightarrow \infty$. 
\end{lemma}

\begin{proof}
It will be convenient to fix a linear operator $B: \R^\ell \twoheadrightarrow \R^{\ell - k}$ such that the linear operator $T := (A, B): \R^\ell \rightarrow \R^\ell$ is an isomorphism. Defining
$$
\mathcal{D}(H) := T^{-1}\left(\prod_{i = 1}^k (c_i, \infty) \times \R^{\ell - k} \right) \cap [-H, H]^\ell,
$$
our task is then to give an asymptotic formula for $\mu(\mathcal{D}(H))$. Changing $c_1$ to $0$ changes the volume by only $O(H^{\ell - 1})$, since the difference is contained in a parallelepiped of side lengths bounded by $O(1), O(H), \dots, O(H)$. Similarly, we may assume that all the $c_i$ are $0$. Let $\mathcal{R}(H)$ denote the region thus obtained, i.e.
$$
\mathcal{R}(H) := T^{-1}\left(\prod_{i = 1}^k (0, \infty) \times \R^{\ell - k} \right) \cap [-H, H]^\ell.
$$
Then $\mathcal{R}(H) = H \mathcal{R}(1)$, so its volume is $H^\ell$ times the volume of $\mathcal{R}(1)$, which is positive since a linear change of coordinates transforms $\mathcal{R}(1)$ into the intersection of $(0, \infty)^k \times \R^{\ell - k}$ with a neighborhood of $0$.
% We have the asymptotic expansion
% $$
% \mu(\mathcal{D}) = \mu\left(T^{-1} \left((0, \infty)^k \times \R^{\ell - k} \right) \cap [-H, H]^\ell \right) + O(H^{\ell - 1})
% $$
% as $H \rightarrow \infty$. Now observe that
% $$
% T^{-1}\left((0, \infty)^k \times \R^{\ell - k} \right) \cap [-H, H]^\ell = H \cdot \left(T^{-1}\left((0, \infty)^k \times \R^{\ell - k} \right) \cap [-1, 1]^\ell\right).
% $$
% Since $T([-1, 1]^\ell)$ is an open set containing $0$, we conclude that $T^{-1}\left((0, \infty)^k \times \R^{\ell - k}\right) \cap [-1, 1]^\ell$ is a non-empty open set and hence has positive Lebesgue measure. Choosing
% $$
% C := \mu\left(T^{-1}\left((0, \infty)^k \times \R^{\ell - k}\right) \cap [-1, 1]^\ell \right) > 0
% $$
% ends the proof of the lemma.
\end{proof}

We use Lemma \ref{lRealVolume} with $\ell = 2n$ and $k$ equal to twice the number of real places of $K$. It is however not obvious that the region $\Omega$ can be given by $k$ linear inequalities, as it is a priori only defined through $2k$ linear inequalities. We will start by establishing a general lemma granting linear independence of the resulting set of linear forms. We write $L_i^{\text{hom}}$ for the homogeneous part of the affine linear polynomial $L_i$.

\begin{lemma} 
\label{lLinearIndependence}
For each real place $\tau \in \Omega_K$, fix a 2-element subset $I_\tau \subseteq \{1, 2, 3, 4\}$. Then $\{\tau \circ L_i^{\textup{hom}} : \tau \textup{ real}, i \in  I_\tau\} \subseteq \Hom_{\Q\textup{-linear}}(K, \R)$ is an $\R$-linearly independent set.
\end{lemma}

\begin{proof}
The four forms $L_i^{\text{hom}}$ are pairwise $K$-linearly independent, so the $L_i^{\text{hom}}$ for $i \in I_\tau$ define an isomorphism $K^2 \rightarrow K^2$. Applying $\otimes_K K_\tau$ and taking the product over all real places $\tau$ gives an
isomorphism $\prod_\tau K_\tau^2 \rightarrow \prod_\tau K_\tau^2$. Composing with the projections onto the various $K_\tau \cong \R$
on the right gives $\R$-linearly independent maps. Pre-composing with the product of two copies of the surjection $K_\R := K \otimes \R \twoheadrightarrow \prod_\tau K_\tau$ gives the $\R$-linear extensions $K_\R^2 \rightarrow \R$ of the maps $\tau \circ L_i^{\text{hom}}$ so the latter are $\R$-linearly independent too.
\end{proof}

Let
$$
c := m^2 \kappa X + 1, \quad \quad d := m^2 \kappa (m^2 \kappa Y + \lambda).
$$
Then observe that $L_i(X, Y) = c + a_id$ for $i \in \{1, 2, 3\}$ and $d = m^2 \kappa L_4(X, Y)$. We shall now use $(K1)$ to rewrite the defining inequalities of $\Omega$ as $k$ linear inequalities, namely $2$ for each real place $\tau$ of $K$, given by two linear forms $M_{\tau, 1}$ and $M_{\tau, 2}$. In order to construct $M_{\tau, 1}$ and $M_{\tau, 2}$, we will now distinguish two cases. 

Firstly, suppose that the real place is $\tau_\ell$ for $\ell \in \{1, \dots, 6\}$. Then we will construct $M_{\tau_\ell, 1}$ and $M_{\tau_\ell, 2}$ via the following lemma.

\begin{lemma}
Let $\ell \in \{1, \dots, 6\}$. Let $I_\ell$ be $\{1, 3\}$ if $\ell \in \{1, 3, 5\}$ and $\{2, 3\}$ if $\ell \in \{2, 4, 6\}$ and define $i_0$ to be the unique element of $\{1, 2, 3\} - I_\ell$. Assume that the sign of $\tau_\ell(L_i(X, Y))$ is precisely the $(i, \ell)$ entry in the table \eqref{eInfMat2} for all $i \in I_\ell$. Then the sign of $\tau_\ell(L_{i_0}(X, Y))$ is also the $(i_0, \ell)$ entry in the table \eqref{eInfMat2}, and moreover $\tau_\ell(L_4(X, Y)) > 0$.
\end{lemma}

\begin{proof}
Appealing to $(K1)$ and $d = m^2 \kappa L_4(X, Y)$, we see that $(IV)$ is equivalent to $\tau_\ell(d) > 0$ for $\ell \in \{1, 2, 3\}$ and $\tau_\ell(d) < 0$ for $\ell \in \{4, 5, 6\}$. If $\{i, j, k\} = \{1, 2, 3\}$, the condition $\tau_\ell(a_i) < \tau_\ell(a_j) < \tau_\ell(a_k)$ implies
\begin{align}
\label{eIneq1}
\tau_\ell(c + a_i d) < \tau_\ell(c + a_j d) < \tau_\ell(c + a_k d)
\end{align}
for $\ell \in \{1, 2, 3\}$, and
\begin{align}
\label{eIneq2}
\tau_\ell(c + a_k d) < \tau_\ell(c + a_j d) < \tau_\ell(c + a_i d)
\end{align}
for $\ell \in \{4, 5, 6\}$. 

Keeping the inequalities \eqref{eIneq1} and \eqref{eIneq2} in mind and the orderings \eqref{eTau12}, \eqref{eTau34}, \eqref{eTau56}, the lemma is now a direct verification.
\end{proof}

If $I_\ell = \{1, 3\}$, we choose $M_{\tau_\ell, 1} := \tau_\ell \circ L_1^{\text{hom}}$ and $M_{\tau_\ell, 2} := \tau_\ell \circ L_3^{\text{hom}}$. If $I_\ell = \{2, 3\}$, we choose $M_{\tau_\ell, 1} := \tau_\ell \circ L_2^{\text{hom}}$ and $M_{\tau_\ell, 2} := \tau_\ell \circ L_3^{\text{hom}}$.

Secondly, suppose that the real place is $\tau$ with $\tau \neq \tau_\ell$ for all $\ell \in \{1, \dots, 6\}$. In this case, we can also rewrite, for each given infinite place, the combined four inequalities from \eqref{eEasyInf1} and \eqref{eEasyInf2} as two inequalities. Indeed, pick $i \in \{1, 2, 3\}$ with $\tau(a_i)$ minimal. Note that $\tau(c + a_id) > 0$ then forces $\tau(c + a_jd) > 0$ for all $j \in \{1, 2, 3\}$ because $\tau(d) > 0$ by \eqref{eEasyInf2}. Then we take $M_{\tau, 1} = \tau \circ L_i^{\text{hom}}$ and $M_{\tau, 2} = \tau \circ L_4^{\text{hom}}$. 

The linear forms $\{M_{\tau, 1} : \tau \in \Omega_K^{\text{real}}\} \cup \{M_{\tau, 2} : \tau \in \Omega_K^{\text{real}}\}$ are linearly independent by Lemma \ref{lLinearIndependence}, so the resulting linear operator $A$ is surjective. By Lemma \ref{lRealVolume}, we conclude that there is some $C > 0$
\begin{align}
\label{eVolumeConclusion}
\mu\left(\Omega \cap [-H, H]^{2n}\right) \sim C H^{2n} \quad \quad \text{as } H \rightarrow \infty.
\end{align}

\subsubsection*{Application of Kai's theorem}
We now apply Theorem \ref{tKai2} to $\Omega$ and the linear polynomials $L_1, L_2, L_3, L_4$. Recall that the linear forms $L_j^{\text{hom}}$ are pairwise linearly independent. Therefore the cokernel condition from Theorem \ref{tKai2} holds. Moreover, the volume assumption in Theorem \ref{tKai2} is satisfied thanks to equation \eqref{eVolumeConclusion}. Observe that the leading constant in Theorem \ref{tKai2} is positive by Lemma \ref{lAdmissible}. 

We now claim that the asymptotic from Theorem \ref{tKai2} remains valid upon replacing the function $\Lambda^{O_K}(L_i(x, y))$ by $\Lambda^{O_K}(L_i(x, y)) \mathbf{1}_{L_i(x, y) \text{ prime element}}$. Indeed, the difference is bounded in absolute value by
\begin{equation}
\label{eFirstBound}
\sum_{j = 1}^4 \sum_{\substack{x, y \in \Z^n \\ (x, y) \in [-N, N]^{2n}}} \mathbf{1}_{(L_j(x, y)) = \mathfrak{p}^m \text{ for some } \mathfrak{p} \in \Omega_K, m \geq 2} \prod_{i = 1}^4 \Lambda^{O_K}(L_i(x, y)).
\end{equation}
Note that we have already fixed $K$, an integral basis for $K$ and the linear forms $L_1, \dots, L_4$. Observe that if $(x, y) \in [-N, N]^{2n}$, then the norm of $L_j(x, y)$ is at most $B_1 N^n$ (with $B_1$ depending only on $K$, the fixed integral basis and $L_1, \dots, L_4$). Then for every $\epsilon > 0$, there exists $B_2(\epsilon) > 0$ (which is allowed to depend on $\epsilon$, $K$, the fixed integral basis of $K$ and $L_1, \dots, L_4$) such that all $(x, y) \in [-N, N]^{2n}$ satisfy $\Lambda^{O_K}(L_i(x, y)) \leq B_2(\epsilon) N^\epsilon$. Therefore the expression in equation \eqref{eFirstBound} is at most
$$
B_2(\epsilon)^4 N^{4\epsilon} \sum_{j = 1}^4 \sum_{\substack{x, y \in \Z^n \\ (x, y) \in [-N, N]^{2n}}} \mathbf{1}_{(L_j(x, y)) = \mathfrak{p}^m \text{ for some } \mathfrak{p} \in \Omega_K, m \geq 2}.
$$
Using once more that the norm of $L_j(x, y)$ is at most $B_1 N^n$ and using that $m \geq 2$, we see that there are at most $B_3 N^{n/2}$ choices for $\mathfrak{p}^m$. We also observe that for every $\epsilon > 0$, there exists $B_4(\epsilon) > 0$ such that for every ideal $I$ of $O_K$ and every $j \in \{1, \dots, 4\}$
$$
\#\{(x, y) \in [-N, N]^{2n} : I = (L_j(x, y))\} \leq B_4(\epsilon) N^{1 + \epsilon}.
$$
Hence equation \eqref{eFirstBound} is at most $\ll_\epsilon N^{\frac{3n}{2} + 5\epsilon}$. This proves the claim, since the leading term in Theorem \ref{tKai2} has magnitude $N^{2n}$.

Furthermore, for any fixed finite set $S \subseteq \Omega_K$, the asymptotic in Theorem \ref{tKai2} also remains valid upon imposing the following additional summation conditions on $(x, y) \in \Z^{2n}$:
\begin{itemize}
\item each prime element $L_i(x, y)$ is not in $S$,
\item the prime elements $L_i(x, y)$ and $L_j(x, y)$ are pairwise coprime.
\end{itemize}
We apply this modified version of Theorem \ref{tKai2} first with $S := T'$, and then we repeatedly apply this with $S$ equal to $T'$ together with the prime elements obtained thus far. Proceeding in this way, we obtain an infinite sequence $(x_j, y_j) \in \Z^{2n}$ such that, upon setting $q_{i, j} := L_i(x_j, y_j)$, we have that
\begin{itemize}
\item each $q_{i, j}$ is a prime element such that $(q_{i, j}) \not \in T'$,
\item for every $(i_1, j_1) \neq (i_2, j_2)$ we have that $(q_{i_1, j_1}) \neq (q_{i_2, j_2})$.
\end{itemize}
Set
$$
c_j := m^2 \kappa x_j + 1, \quad \quad d_j := m^2 \kappa (m^2 \kappa y_j + \lambda)
$$
and
$$
t_j := \kappa q_{1, j} q_{2, j} q_{3, j} q_{4, j} = m^{-2} d_j (c_j + a_1 d_j) (c_j + a_2 d_j) (c_j + a_3 d_j).
$$
We claim that $t_j$ satisfies $(P1)$, $(P2)$, $(P3)$. Moreover, we claim that $t_j$ satisfies $(P4)$ for all but finitely many quadruples. Once we establish the claim, we find in particular at least one suitable $t_j$, proving Theorem \ref{tSuitable}. Thus it remains to prove the claim. In order to lighten up the notation, we fix $j$ for now and write $t$, $c$, $d$, $q_1$, $q_2$, $q_3$, $q_4$ for $t_j$, $c_j$, $d_j$, $q_{1, j}$, $q_{2, j}$, $q_{3, j}$, $q_{4, j}$.

Let us first check $(P1)$. By construction of the $L_i$, we see that 
$$
q_1 \equiv q_2 \equiv q_3 \equiv 1 \bmod 8N \kappa.
$$
Moreover, by the definition of $\lambda$ and $L_4$, we have $q_4 \kappa \equiv \lambda \kappa \equiv 1 \bmod 8N$, hence $t \equiv 1 \bmod 8N$, so by Hensel’s lemma, $t \in K_v^{\ast 2}$ for all finite $v \in T$. We now check that $t$ is totally positive. If $\tau$ is a real place other than $\tau_1, \dots, \tau_6$, then $\tau(\kappa q_4) > 0$ by $(II)$, and $\tau(q_1), \tau(q_2), \tau(q_3) > 0$ by $(I)$; hence $\tau(t) > 0$. If instead $\tau = \tau_\ell$ for some $\ell \in \{1, \dots, 6\}$, then $\tau(\kappa q_4)$ is positive for $\ell \in \{1, 2, 3\}$ and negative for $\ell \in \{4, 5, 6\}$ by $(IV)$ and $(K1)$; by inspecting the table \eqref{eInfMat2}, this matches exactly the sign of $\tau(q_1 q_2 q_3)$, so $\tau(t) > 0$. Hence $t$ is totally positive. This establishes $(P1)$. Finally, the identity $t = \kappa q_1 q_2 q_3 q_4$ gives $(P2)$.

We now check $(P3)$. We have shown just above that
$$
q_1 \equiv q_2 \equiv q_3 \equiv 1 \bmod 8N \kappa
$$
and that $q_1, q_2, q_3$ are positive at all real embeddings outside $\tau_1, \dots, \tau_6$, so $q_1, q_2, q_3 \in K_v^{\ast 2}$ for all $v \in T' - \{\tau_1, \dots, \tau_6\}$. Now fix $\boldsymbol{\pi}$ and $z_i$ as in $(K2)$. Then we have
$$
\prod_{v \in T'} (z_i, q_j)_v = \prod_{v \mid \infty} (z_i, q_j)_v = \prod_{v \in \{\tau_1, \tau_2, \tau_3, \tau_4, \tau_5, \tau_6\}} (z_i, q_j)_v
$$
for $i \in \{1, \dots, 12\}$ and $j \in \{1, 2, 3\}$. By $(K2)$, the signs of $z_i$ at $\tau_1, \dots, \tau_6$ are given by the table \eqref{eInfMat1}. Then $(P3)$ follows from a direct computation using the tables \eqref{eInfMat1} and \eqref{eInfMat2}.

It remains to check $(P4)$. Because $t$ is of the shape $m^{-2} d (c + a_1 d) (c + a_2 d) (c + a_3 d)$, one directly checks that
$$
E^{-t}: m^{-2} d (c + a_1 d) (c + a_2 d) (c + a_3 d) y^2 = (x + a_1) (x + a_2) (x + a_3)
$$
has the rational point $(x, y) = (c/d, m/d^2)$, which is also readily seen to be non-torsion for all but finitely many quadruples $(q_{1, j}, q_{2, j}, q_{3, j}, q_{4, j})$ by Lemma \ref{lTorsion}.
\end{proof}

\subsection{Construction of an auxiliary twist}
\label{ssFinal}
It remains to prove that an auxiliary twist $\kappa$ exists, i.e.~some $\kappa$ satisfying $(K1)$ and $(K2)$. It is possible to give a probabilistic proof of this statement, in the spirit of \cite{SD}. However, for the sake of brevity, we have opted for a direct construction. Before we start with the proof, let us state an important technical result. This lemma classifies precisely in which situation one is unable to reduce the dimension by $2$ in the context of Lemma \ref{lSelmerChange}.

\begin{lemma}
\label{lReduction}
Let $V \subseteq H^1(G_K, \mathbb{F}_2^2)$ be a subspace with $\dim_{\FF_2} V \geq 2$. Assume that there exists a finite set of places $S \subseteq \Omega_K$ such that for all $(\chi_1, \chi_2), (\chi_1', \chi_2') \in V$ and all finite places $\mathfrak{p} \not \in S$
$$
\det
\begin{pmatrix}
\chi_1(\Frob_\mathfrak{p}) & \chi_2(\Frob_\mathfrak{p}) \\
\chi_1'(\Frob_\mathfrak{p}) & \chi_2'(\Frob_\mathfrak{p})
\end{pmatrix}
= 0.
$$
Then, writing $T_1, T_2: H^1(G_K, \mathbb{F}_2^2) \rightarrow H^1(G_K, \mathbb{F}_2)$ for the natural projection maps, we have $T_1(V) = 0$, $T_2(V) = 0$ or $(T_1 + T_2)(V) = 0$.
\end{lemma}

\begin{proof}
Since $\ker(T_1)$, $\ker(T_2)$ and $\ker(T_1 + T_2)$ have pairwise trivial intersection, we reduce to the case that $\dim_{\FF_2} V = 2$, so $(\chi_1, \chi_2)$ and $(\chi_1', \chi_2')$ form a basis of $V$.

The Chebotarev density theorem and our assumption imply that the homomorphism $\varphi: G_K \rightarrow M_2(\mathbb{F}_2)$ given by
$$
\sigma \mapsto 
\begin{pmatrix}
\chi_1(\sigma) & \chi_2(\sigma) \\
\chi_1'(\sigma) & \chi_2'(\sigma)
\end{pmatrix}
$$
satisfies $\det \circ \varphi = 0$. Write $W := \varphi(G_K)$. We observe that $(M_2(\mathbb{F}_2), \det)$ is a weakly metabolic quadratic space over $\mathbb{F}_2$ and that $W$ is an isotropic subspace. Since there are exactly six maximal isotropic subspaces by \cite[Proposition 2.6(b)]{PR}, it follows that $W$ is included in one of the following six subspaces of $M_2(\mathbb{F}_2)$
$$
\left\{ 
\begin{pmatrix}
* & * \\
0 & 0
\end{pmatrix}
,
\begin{pmatrix}
0 & 0 \\
* & *
\end{pmatrix}
,
\begin{pmatrix}
* & 0 \\
* & 0
\end{pmatrix}
,
\begin{pmatrix}
0 & * \\
0 & *
\end{pmatrix}
,
\begin{pmatrix}
a & a \\
b & b
\end{pmatrix}
,
\begin{pmatrix}
a & b \\
a & b
\end{pmatrix}
\right\}
.
$$
Since $(\chi_1, \chi_2)$ and $(\chi_1', \chi_2')$ are linearly independent, it is easy to eliminate the first, second and sixth possibility. The remaining three options correspond exactly to the final conclusion of the lemma.
\end{proof}

We will prove the existence of an auxiliary twist in four parts. In the first part, we construct distinct prime ideals $\mathfrak{p}_1, \dots, \mathfrak{p}_{s_1} \not \in T$ and corresponding local uniformizers $\pi_1, \dots, \pi_{s_1}$ for some integer $s_1 \geq 0$ such that
$$
\mathrm{Sel}_{\mathcal{L}_{s_1, (\pi_1, \dots, \pi_{s_1})}}(G_K, E[2]) = 1. 
$$
In the second part, we enlarge this list of prime ideals to $\mathfrak{p}_1, \dots, \mathfrak{p}_{s_2} \not \in T$ and corresponding local uniformizers $\pi_1, \dots, \pi_{s_2}$ for some integer $s_2 \geq s_1$ such that certain auxiliary spaces vanish while maintaining the property $\mathrm{Sel}_{\mathcal{L}_{s_2, (\pi_1, \dots, \pi_{s_2})}}(G_K, E[2]) = 1$.

In the third and fourth part, we introduce $16$ more prime ideals (in particular $s := s_2 + 16$). We will first construct six convenient prime elements $p_{s_2 + 1}, \dots, p_{s_2 + 6}$ with prescribed local behavior and local uniformizers. As we are introducing these prime elements, our Selmer groups may grow and we use six additional prime ideals $\mathfrak{p}_{s_2 + 7}, \dots, \mathfrak{p}_{s_2 + 12}$ to remove any additional classes introduced in this process. It is at this stage that we heavily rely on the auxiliary vanishing from the second part.

Finally, we introduce three more prime elements $p_{s_2 + 13}$, $p_{s_2 + 14}$ and $p_{s_2 + 15}$, which combined appropriately with $p_{s_2 + 1}, \dots, p_{s_2 + 6}$, introduce Selmer classes satisfying the sign conditions from the table \eqref{eInfMat1}. The task of $\mathfrak{p}_{s_2 + 16}$ is to ensure the sign conditions and ramification conditions of $(K1)$ and to principalize the product of all the prime ideals. 

\subsubsection{Reducing ranks}
Throughout this section, $\boldsymbol{\pi}$ denotes a vector of uniformizers (the precise entries of this vector will be clear from context). We start by determining the parity of $\dim_{\FF_2} \mathrm{Sel}_{\mathcal{L}_{i, \boldsymbol{\pi}}}(G_K, E[2])$.

\begin{lemma}
\label{lParity}
We have for all $i \geq 0$, for all distinct prime ideals $\mathfrak{p}_1, \dots, \mathfrak{p}_i \not \in T$ and all local uniformizers $\pi_1, \dots, \pi_i$
$$
\dim_{\FF_2} \mathrm{Sel}_{\mathcal{L}_{i, \boldsymbol{\pi}}}(G_K, E[2]) \equiv 0 \bmod 2.
$$
\end{lemma}

\begin{proof}
For $i = 0$, we observe that $\mathrm{Sel}_{\mathcal{L}_{0, \boldsymbol{\pi}}}(G_K, E[2]) = \mathrm{Sel}^2(E/K)$ by Lemma \ref{lInitialState}. Hence the case $i = 0$ follows from Corollary \ref{cCes} and Lemma \ref{lStartingCurve}$(c)$. The general case is then a consequence of Lemma \ref{lSelmerChange} and the case $i = 0$.
\end{proof}

\begin{lemma}
\label{lCheb}
For every two non-trivial elements $y_1, y_2 \in K^\ast/K^{\ast 2}$, there exist infinitely many prime ideals $\mathfrak{q}$ such that $y_1$ and $y_2$ are non-trivial in $K_{\mathfrak{q}}^\ast/K_{\mathfrak{q}}^{\ast 2}$. 
\end{lemma}

\begin{proof}
This follows from the Chebotarev density theorem if $y_1$ and $y_2$ are linearly independent. If $y_1$ and $y_2$ are not linearly independent, then we must have $y_1 = y_2$ and moreover $y_1$ is non-trivial. In this case the lemma follows once more from the Chebotarev density theorem.
\end{proof}

Our next lemma corresponds precisely to the first part of the argument outlined above. 

\begin{lemma}
\label{lAux1}
There exists an integer $s_1 \geq 0$, distinct prime ideals $\mathfrak{p}_1, \dots, \mathfrak{p}_{s_1} \not \in T$ and local uniformizers $\pi_1, \dots, \pi_{s_1}$ with 
$$
\mathrm{Sel}_{\mathcal{L}_{s_1, \boldsymbol{\pi}}}(G_K, E[2]) = 1. 
$$
\end{lemma}

\begin{proof}
In order to prove the lemma, we claim that for all distinct primes $\mathfrak{p}_1, \dots, \mathfrak{p}_i \not \in T$ and for all local uniformizers $\pi_1, \dots, \pi_i$ at least one of the following three statements hold:
\begin{enumerate}
\item[$(S1)$] we have $\mathrm{Sel}_{\mathcal{L}_{i, \boldsymbol{\pi}}}(G_K, E[2]) = 1$,
\item[$(S2)$] there exist a prime ideal $\mathfrak{p}_{i + 1}$ and a local uniformizer $\pi_{i + 1}$ such that
$$
\dim_{\FF_2} \mathrm{Sel}_{\mathcal{L}_{i + 1, \boldsymbol{\pi}}}(G_K, E[2]) = - 2 + \dim_{\FF_2} \mathrm{Sel}_{\mathcal{L}_{i, \boldsymbol{\pi}}}(G_K, E[2]),
$$
\item[$(S3)$] there exist prime ideals $\mathfrak{p}_{i + 1}, \mathfrak{p}_{i + 2}$ and local uniformizers $\pi_{i + 1}, \pi_{i + 2}$ such that
$$
\dim_{\FF_2} \mathrm{Sel}_{\mathcal{L}_{i + 2, \boldsymbol{\pi}}}(G_K, E[2]) = - 2 + \dim_{\FF_2} \mathrm{Sel}_{\mathcal{L}_{i, \boldsymbol{\pi}}}(G_K, E[2]).
$$
\end{enumerate}
The claim visibly implies Lemma \ref{lAux1}, so it remains to prove the claim.

Let us now suppose that $(S1)$ fails. Thus we have $\dim_{\FF_2} \mathrm{Sel}_{\mathcal{L}_{i, \boldsymbol{\pi}}}(G_K, E[2]) \geq 2$ by Lemma \ref{lParity}. If the space $V := \mathrm{Sel}_{\mathcal{L}_{i, \boldsymbol{\pi}}}(G_K, E[2])$ does not satisfy the hypothesis of Lemma \ref{lReduction}, then we can find a prime ideal $\mathfrak{p}_{i + 1} \not \in T \cup \{\mathfrak{p}_1, \dots, \mathfrak{p}_i\}$ and two elements $(\chi_1, \chi_2) \in V$ and $(\chi_1', \chi_2') \in V$ such that
\begin{equation}
\label{eDet1}
\det
\begin{pmatrix}
\chi_1(\Frob_{\mathfrak{p}_{i + 1}}) & \chi_2(\Frob_{\mathfrak{p}_{i + 1}}) \\
\chi_1'(\Frob_{\mathfrak{p}_{i + 1}}) & \chi_2'(\Frob_{\mathfrak{p}_{i + 1}})
\end{pmatrix}
= 1.
\end{equation}
Fix an arbitrary choice for the local uniformizer $\pi_{i + 1}$. By equation \eqref{eDet1}, the restrictions $\res_{\mathfrak{p}_{i + 1}} (\chi_1, \chi_2)$ and $\res_{\mathfrak{p}_{i + 1}} (\chi_1', \chi_2')$ are linearly independent, hence Remark \ref{rSelmerChangeUpgrade}$(c)$ shows that we must be in the second case of Lemma \ref{lSelmerChange}. Thus $(S2)$ holds in this case.

We henceforth assume that $V := \mathrm{Sel}_{\mathcal{L}_{i, \boldsymbol{\pi}}}(G_K, E[2])$ does satisfy the hypothesis of Lemma \ref{lReduction}. Therefore, recalling that $T_1, T_2$ are the projection maps from Lemma \ref{lReduction}, it follows from Lemma \ref{lReduction} that 
$$
T_1(\mathrm{Sel}_{\mathcal{L}_{i, \boldsymbol{\pi}}}(G_K, E[2])) = 0 \ \text{ or } \ T_2(\mathrm{Sel}_{\mathcal{L}_{i, \boldsymbol{\pi}}}(G_K, E[2])) = 0 \ \text{ or } \ (T_1 + T_2)(\mathrm{Sel}_{\mathcal{L}_{i, \boldsymbol{\pi}}}(G_K, E[2])) = 0.
$$
Our goal is now to show that $(S3)$ is true. We proceed with the case $T_1(\mathrm{Sel}_{\mathcal{L}_{i, \boldsymbol{\pi}}}(G_K, E[2])) = 0$, the other cases are similar.

Recalling that $\dim_{\FF_2} \mathrm{Sel}_{\mathcal{L}_{i, \boldsymbol{\pi}}}(G_K, E[2]) \geq 2$, we may therefore fix some non-trivial element $(1, x) \in \mathrm{Sel}_{\mathcal{L}_{i, \boldsymbol{\pi}}}(G_K, E[2])$. We apply Lemma \ref{lCheb} with $y_1 = x$ and $y_2 = \alpha \beta$; note that $x$ is not a square by assumption and that $\alpha \beta$ is not a square by Lemma \ref{lStartingCurve}$(a)$. Therefore there exists a prime ideal $\mathfrak{p}_{i + 1} \not \in T \cup \{\mathfrak{p}_1, \dots, \mathfrak{p}_i\}$ such that $x$ and $\alpha \beta$ are both non-trivial elements in $K_{\mathfrak{p}_{i + 1}}^\ast/K_{\mathfrak{p}_{i + 1}}^{\ast 2}$. Now we claim that
\begin{equation}
\label{eSelRed1}
\dim_{\FF_2} \res_{\mathfrak{p}_{i + 1}} \, \mathrm{Sel}_{\mathcal{L}_{i, \boldsymbol{\pi}}}(G_K, E[2]) = 1.
\end{equation}
Indeed, the above dimension is at least $1$, as the element $\res_{\mathfrak{p}_{i + 1}} \, (1, x)$ is non-trivial by the choice of $\mathfrak{p}_{i + 1}$. Moreover, the above dimension is also at most $1$, as $\res_{\mathfrak{p}_{i + 1}} \, \mathrm{Sel}_{\mathcal{L}_{i, \boldsymbol{\pi}}}(G_K, E[2])$ is contained in the subspace $\ker(T_1) \cap H^1_{\text{nr}}(G_{K_{\mathfrak{p}_{i + 1}}}, E[2])$ (which is the $1$-dimensional space $\{0\} \times H^1_{\text{nr}}(G_{K_{\mathfrak{p}_{i + 1}}}, \mathbb{F}_2)$ using our identification $E[2] \cong \mathbb{F}_2^2$ as Galois modules), thus proving the claim \eqref{eSelRed1}. 

We now pick any local uniformizer $\pi_{i + 1}$. It follows from equation \eqref{eSelRed1} that, irregardless of the choice of $\pi_{i + 1}$, we are in the third case of Lemma \ref{lSelmerChange}, and hence
\begin{equation}
\label{e1xGone2}
\dim_{\FF_2} \mathrm{Sel}_{\mathcal{L}_{i + 1, \boldsymbol{\pi}}}(G_K, E[2]) = \dim_{\FF_2} \mathrm{Sel}_{\mathcal{L}_{i, \boldsymbol{\pi}}}(G_K, E[2]).
\end{equation}
We claim that there exists $(x_1, x_2) \in \mathrm{Sel}_{\mathcal{L}_{i + 1, \boldsymbol{\pi}}}(G_K, E[2]) - \ker(T_1)$. Recall that we picked $x$ to be a non-trivial element in $K_{\mathfrak{p}_{i + 1}}^\ast/K_{\mathfrak{p}_{i + 1}}^{\ast 2}$ and that $(1, x) \in \mathrm{Sel}_{\mathcal{L}_{i, \boldsymbol{\pi}}}(G_K, E[2])$. Hence $\res_{\mathfrak{p}_{i + 1}} \, (1, x)$ is a non-trivial element of $H^1_{\text{nr}}(G_{K_{\mathfrak{p}_{i + 1}}}, E[2])$. We deduce that $\res_{\mathfrak{p}_{i + 1}} \, (1, x) \not \in \mathcal{L}_{i + 1, \boldsymbol{\pi}, \mathfrak{p}_{i + 1}}$ by Lemma \ref{lEasy}$(d)$ and therefore
\begin{equation}
\label{e1xGone}
(1, x) \not \in \mathrm{Sel}_{\mathcal{L}_{i + 1, \boldsymbol{\pi}}}(G_K, E[2]).
\end{equation}
By equations \eqref{e1xGone2}, \eqref{e1xGone} and the fact that $(1, x) \in \mathrm{Sel}_{\mathcal{L}_{i, \boldsymbol{\pi}}}(G_K, E[2])$, we conclude that $\mathrm{Sel}_{\mathcal{L}_{i + 1, \boldsymbol{\pi}}}(G_K, E[2]) - \mathrm{Sel}_{\mathcal{L}_{i, \boldsymbol{\pi}}}(G_K, E[2])$ is non-empty. Now we fix any element $(x_1, x_2) \in \mathrm{Sel}_{\mathcal{L}_{i + 1, \boldsymbol{\pi}}}(G_K, E[2]) - \mathrm{Sel}_{\mathcal{L}_{i, \boldsymbol{\pi}}}(G_K, E[2])$, and we will show that $x_1$ is not locally trivial at $\mathfrak{p}_{i + 1}$; this implies that $x_1$ is also not globally trivial and thus $(x_1, x_2) \not \in \ker(T_1)$ as claimed.

Observe that $(x_1, x_2)$ must certainly be non-trivial locally at $\mathfrak{p}_{i + 1}$, as otherwise $(x_1, x_2)$ is in $\mathrm{Sel}_{\mathcal{L}_{i, \boldsymbol{\pi}}}(G_K, E[2])$ contrary to our assumption. Inspecting the three non-trivial elements of $\mathcal{L}_{i + 1, \boldsymbol{\pi}, \mathfrak{p}_{i + 1}}$ in equation \eqref{eSelmerStructure}, we see that $\mathrm{res}_{\mathfrak{p}_{i + 1}}(x_1)$ has odd valuation in two cases (and is hence non-trivial) and that $\mathrm{res}_{\mathfrak{p}_{i + 1}}(x_1) = \mathrm{res}_{\mathfrak{p}_{i + 1}}(\alpha \beta)$ in the third case (which is non-trivial by our choice of $\mathfrak{p}_{i + 1}$). This proves the claimed existence of $(x_1, x_2)$, and we fix any such choice.

We also clearly have
$$
\ker\left( \mathrm{Sel}_{\mathcal{L}_{i, \boldsymbol{\pi}}}(G_K, E[2]) \xrightarrow{\res_{\mathfrak{p}_{i + 1}}} H^1(G_{K_{\mathfrak{p}_{i + 1}}}, E[2]) \right) \subseteq \mathrm{Sel}_{\mathcal{L}_{i + 1, \boldsymbol{\pi}}}(G_K, E[2]).
$$
By the inequality $\dim_{\mathbb{F}_2} \mathrm{Sel}_{\mathcal{L}_{i, \boldsymbol{\pi}}}(G_K, E[2]) \geq 2$ and by equation \eqref{eSelRed1}, there exists a non-trivial element in the above kernel. By the assumption $T_1(\mathrm{Sel}_{\mathcal{L}_{i, \boldsymbol{\pi}}}(G_K, E[2])) = 0$, we see that this element must be of the shape $(1, x')$. 

By Lemma \ref{lCheb}, there exists a prime ideal $\mathfrak{p}_{i + 2} \not \in T \cup \{\mathfrak{p}_1, \dots, \mathfrak{p}_{i + 1}\}$ such that $x_1$ and $x'$ are not locally trivial at $\mathfrak{p}_{i + 2}$. Hence the restrictions of $(1, x'), (x_1, x_2) \in \mathrm{Sel}_{\mathcal{L}_{i + 1, \boldsymbol{\pi}}}(G_K, E[2])$ at $\mathfrak{p}_{i + 2}$ must be linearly independent. By Remark \ref{rSelmerChangeUpgrade}$(c)$ we conclude that we are necessarily in the second case of Lemma \ref{lSelmerChange} for any choice of $\pi_{i + 2}$. We have thus shown that $(S3)$ holds.
\end{proof}

\subsubsection{Vanishing of auxiliary spaces}
We have now completed the first part of the construction of the auxiliary twist. The second part also imposes vanishing of some auxiliary spaces. To state it formally, we introduce for distinct primes $\mathfrak{p}_1, \dots, \mathfrak{p}_r \not \in T$ and local uniformizers $\pi_1, \dots, \pi_r$ the following spaces
$$
\mathcal{D}_{r, 1} := \left\{x \in K^\ast/K^{\ast 2} : 
\begin{array}{c} 
\res_{\mathfrak{p}_j}(x, x) \in \mathcal{L}_{r, \boldsymbol{\pi}, \mathfrak{p}_j} \text{ for all } 1 \leq j \leq r, \\
v(x) \equiv 0 \bmod 2 \text{ for all } v \not \in T \cup \{\mathfrak{p}_1, \dots, \mathfrak{p}_r\}
\end{array}
\right\}
$$
and
$$
\mathcal{D}_{r, 2} := \left\{x \in K^\ast/K^{\ast 2} : 
\begin{array}{c} 
\res_{\mathfrak{p}_j}(x, 1) \in \mathcal{L}_{r, \boldsymbol{\pi}, \mathfrak{p}_j} \text{ for all } 1 \leq j \leq r, \\
v(x) \equiv 0 \bmod 2 \text{ for all } v \not \in T \cup \{\mathfrak{p}_1, \dots, \mathfrak{p}_r\}
\end{array}
\right\}
$$
and
$$
\mathcal{D}_{r, 3} := \left\{x \in K^\ast/K^{\ast 2} : 
\begin{array}{c} 
\res_{\mathfrak{p}_j}(1, x) \in \mathcal{L}_{r, \boldsymbol{\pi}, \mathfrak{p}_j} \text{ for all } 1 \leq j \leq r, \\
v(x) \equiv 0 \bmod 2 \text{ for all } v \not \in T \cup \{\mathfrak{p}_1, \dots, \mathfrak{p}_r\}
\end{array}
\right\}.
$$

\begin{lemma}
\label{lAux2}
There exists an integer $s_2 \geq 0$, distinct prime ideals $\mathfrak{p}_1, \dots, \mathfrak{p}_{s_2} \not \in T$ and local uniformizers $\pi_1, \dots, \pi_{s_2}$ with 
\begin{align}
\label{eDiagonalGone}
\mathrm{Sel}_{\mathcal{L}_{s_2, \boldsymbol{\pi}}}(G_K, E[2]) = 1, \quad \quad \mathcal{D}_{s_2 ,i} = 1 \quad \textup{ for all } i \in \{1, 2, 3\}.
\end{align}
\end{lemma}

\begin{proof}
Fix an integer $s_1 \geq 0$, distinct prime ideals $\mathfrak{p}_1, \dots, \mathfrak{p}_{s_1} \not \in T$ and local uniformizers $\pi_1, \dots, \pi_{s_1}$ as guaranteed by Lemma \ref{lAux1}. We will now find an integer $s_2 \geq s_1$ and enlargements $\mathfrak{p}_1, \dots, \mathfrak{p}_{s_2}$ and $\pi_1, \dots, \pi_{s_2}$ of our previous list of prime ideals and local uniformizers such that equation \eqref{eDiagonalGone} holds.

As a first step, we add four distinct prime ideals $\mathfrak{p}_{s_1 + 1}, \dots, \mathfrak{p}_{s_1 + 4} \not \in T \cup \{\mathfrak{p}_1, \dots, \mathfrak{p}_s\}$ such that every non-trivial $x \in \langle -1, \alpha, \beta, \gamma \rangle$ reduces to a non-square unit modulo at least one $\mathfrak{p}_{s_1 + j}$. By our choice of $T$, the elements $-1, \alpha, \beta, \gamma$ have even valuation at all places $v \not \in T$. Hence 
\begin{equation}
\label{e1abc}
\langle -1, \alpha, \beta, \gamma \rangle \cap \mathcal{D}_{s_1 + 4, i} = 1
\end{equation}
by Lemma \ref{lEasy}$(d)$. We recall the definition of $A_{\mathfrak{p}_{s_1 + j}}$ in equation \eqref{eRelaxedA}, and we choose the local uniformizer $\pi_{s_1 + j}$ each time in such a way that 
\begin{equation}
\label{eUniPick}
\mathcal{L}_{s_1 + j, \boldsymbol{\pi}, \mathfrak{p}_{s_1 + j}} \neq A_{\mathfrak{p}_{s_1 + j}}
\end{equation}
for $j \in \{1, \dots, 4\}$ (this is indeed possible; if we let $\pi_{s_1 + j}$ and $\pi_{s_1 + j}'$ denote the two different local uniformizers in $K_{\mathfrak{p}_{s_ 1 + j}}^\ast/K_{\mathfrak{p}_{s_1 + j}}^{\ast 2}$ and if we denote by $\boldsymbol{\pi}$ and $\boldsymbol{\pi'}$ the resulting vectors, then an inspection of equation \eqref{eSelmerStructure} yields $\mathcal{L}_{s_1 + j, \boldsymbol{\pi}, \mathfrak{p}_{s_1 + j}} \neq \mathcal{L}_{s_1 + j, \boldsymbol{\pi}', \mathfrak{p}_{s_1 + j}}$). Since the restriction of $\mathrm{Sel}_{\mathcal{L}_{s_1, \boldsymbol{\pi}}}(G_K, E[2]) = 1$ at $\mathfrak{p}_{s_1 + 1}$ is trivial, we are not in the second case of Lemma \ref{lSelmerChange}, and the first case of Lemma \ref{lSelmerChange} is ruled out by equation \eqref{eUniPick}. Therefore we must be in the third case of Lemma \ref{lSelmerChange}. Hence $\mathrm{Sel}_{\mathcal{L}_{s_1 + 1, \boldsymbol{\pi}}}(G_K, E[2]) = 1$ by Lemma \ref{lSelmerChange}. Repeating this argument three more times, we deduce that
\begin{align}
\label{eRayClass}
\mathrm{Sel}_{\mathcal{L}_{s_1 + 4, \boldsymbol{\pi}}}(G_K, E[2]) = 1.
\end{align}
% We may assume that 
% $$
% \dim_{\FF_2} \mathcal{D}_{s_1 + 4, 1} + \dim_{\FF_2} \mathcal{D}_{s_1 + 4, 2} + \dim_{\FF_2} \mathcal{D}_{s_1 + 4, 3} > 0,
% $$
% as otherwise we can take $s_2 := s_1 + 4$. 
We claim that for all $j \geq 0$ and for all distinct prime ideals $\mathfrak{p}_{s_1 + 5}, \dots, \mathfrak{p}_{s_1 + 4 + j} \not \in T \cup \{\mathfrak{p}_1, \dots, \mathfrak{p}_{s_1 + 4}\}$ and for all local uniformizers $\pi_{s_1 + 5}, \dots, \pi_{s_1 + 4 + j}$ satisfying 
\begin{gather}
\mathrm{Sel}_{\mathcal{L}_{s_1 + 4 + j, \boldsymbol{\pi}}}(G_K, E[2]) = 1 \quad \text{ and } \quad \sum_{i = 1}^3 \dim_{\FF_2} \mathcal{D}_{s_1 + 4 + j, i} > 0, \nonumber \\
\mathcal{D}_{s_1 + 4 + j, i} \subseteq \mathcal{D}_{s_1 + 4, i} \quad \textup{ for all } i \in \{1, 2, 3\}, \label{ePreImp}
\end{gather}
there exists a prime ideal $\mathfrak{p}_{s_1 + 5 + j}$ and a local uniformizer $\pi_{s_1 + 5 + j}$ such that 
\begin{gather}
\mathrm{Sel}_{\mathcal{L}_{s_1 + 5 + j, \boldsymbol{\pi}}}(G_K, E[2]) = 1 \quad \text{ and } \quad \sum_{i = 1}^3 \dim_{\FF_2} \mathcal{D}_{s_1 + 5 + j, i} < \sum_{i = 1}^3 \dim_{\FF_2} \mathcal{D}_{s_1 + 4 + j, i}, \nonumber \\
\mathcal{D}_{s_1 + 5 + j, i} \subseteq \mathcal{D}_{s_1 + 4, i} \quad \textup{ for all } i \in \{1, 2, 3\}. \label{ePostImp}
\end{gather}
Note that this claim suffices to prove the lemma; if equation \eqref{ePreImp} fails for $j = 0$, then we must have
$$
\sum_{i = 1}^3 \dim_{\FF_2} \mathcal{D}_{s_1 + 4 + j, i} = 0
$$
by equation \eqref{eRayClass} and hence we may pick $s_2 := s_1 + 4$, while if equation \eqref{ePreImp} holds for $j = 0$, we may iteratively apply this claim to obtain the lemma.

We will now prove the claim. Let $j \geq 0$, let $\mathfrak{p}_{s_1 + 5}, \dots, \mathfrak{p}_{s_1 + 4 + j} \not \in T \cup \{\mathfrak{p}_1, \dots, \mathfrak{p}_{s_1 + 4}\}$ be distinct and let $\pi_{s_1 + 5}, \dots, \pi_{s_1 + 4 + j}$ be local uniformizers satisfying equation \eqref{ePreImp}. Fix some $i$ and some non-trivial $x \in \mathcal{D}_{s_1 + 4 + j, i}$. Since $x \not \in \langle -1, \alpha, \beta, \gamma \rangle$ by \eqref{ePreImp} and \eqref{e1abc}, it follows from Lemma \ref{lStartingCurve}$(a)$ that $\langle -1, \alpha, \beta, \gamma, x \rangle$ spans a $5$-dimensional subspace of $K^\ast/K^{\ast 2}$. Therefore there exists a prime ideal $\mathfrak{p}_{s_1 + 5 + j} \not \in T \cup \{\mathfrak{p}_1, \dots, \mathfrak{p}_{s_1 + 4 + j}\}$ such that
\begin{equation}
\label{eResabc}
\res_{\mathfrak{p}_{s_1 + 5 + j}}(x) \neq 1, \quad \res_{\mathfrak{p}_{s_1 + 5 + j}}(\alpha \beta) \neq 1, \quad \res_{\mathfrak{p}_{s_1 + 5 + j}}(-\alpha \gamma) \neq 1, \quad \res_{\mathfrak{p}_{s_1 + 5 + j}}(\beta \gamma) \neq 1.
\end{equation}
We recall the definition of $A_{\mathfrak{p}_{s_1 + 5 + j}}$ in equation \eqref{eRelaxedA}, and we choose the local uniformizer $\pi_{s_1 + 5 + j}$ in such a way that 
$$
\mathcal{L}_{s_1 + 5 + j, \boldsymbol{\pi}, \mathfrak{p}_{s_1 + 5 + j}} \neq A_{\mathfrak{p}_{s_1 + 5 + j}}.
$$
Since $\mathrm{Sel}_{\mathcal{L}_{s_1 + 4 + j, \boldsymbol{\pi}}}(G_K, E[2]) = 1$ by equation \eqref{ePreImp}, we may use the argument after equation \eqref{eUniPick} to show that we are in the third case of Lemma \ref{lSelmerChange}. This guarantees that 
\begin{equation}
\label{eSelReminder}
\mathrm{Sel}_{\mathcal{L}_{s_1 + 5 + j, \boldsymbol{\pi}}}(G_K, E[2]) = 1. 
\end{equation}
Observe that 
\begin{equation}
\label{exReminder}
x \in \mathcal{D}_{s_1 + 4 + j, i} - \mathcal{D}_{s_1 + 5 + j, i};
\end{equation}
indeed, $x$ has even valuation at $\mathfrak{p}_{s_1 + 5 + j}$ and is non-trivial by equation \eqref{eResabc}, so our observation follows from Lemma \ref{lEasy}$(d)$. 

We will now show that $\mathcal{D}_{s_1 + 5 + j, i'} \subseteq \mathcal{D}_{s_1 + 4 + j, i'}$ for each $i'$. To this end, we claim that by our choice of $v := \mathfrak{p}_{s_1 + 5 + j}$
\begin{align}
&\mathcal{L}_{s_1 + 5 + j, \boldsymbol{\pi}, \mathfrak{p}_{s_1 + 5 + j}} \cap \{(y, y) : y \in K_v^\ast/K_v^{\ast 2}\} = 1, \label{eIntersect1} \\
&\mathcal{L}_{s_1 + 5 + j, \boldsymbol{\pi}, \mathfrak{p}_{s_1 + 5 + j}} \cap \{(y, 1) : y \in K_v^\ast/K_v^{\ast 2}\} = 1, \label{eIntersect2} \\
&\mathcal{L}_{s_1 + 5 + j, \boldsymbol{\pi}, \mathfrak{p}_{s_1 + 5 + j}} \cap \{(1, y) : y \in K_v^\ast/K_v^{\ast 2}\} = 1. \label{eIntersect3}
\end{align}
Indeed, equation \eqref{eIntersect1} follows from directly comparing the three non-trivial elements in $\{(y, y) : y \in K_v^\ast/K_v^{\ast 2}\}$ with $\mathcal{L}_{s_1 + 5 + j, \boldsymbol{\pi}, \mathfrak{p}_{s_1 + 5 + j}}$ and using that $\res_{\mathfrak{p}_{s_1 + 5 + j}}(\beta \gamma) \neq 1$ from equation \eqref{eResabc}. The other claims \eqref{eIntersect2} and \eqref{eIntersect3} are proven similarly.

By equations \eqref{eIntersect1}, \eqref{eIntersect2} and \eqref{eIntersect3}, every element in $\mathcal{D}_{s_1 + 5 + j, i'}$ must be trivial at $\mathfrak{p}_{s_1 + 5 + j}$. This gives the inclusion $\mathcal{D}_{s_1 + 5 + j, i'} \subseteq \mathcal{D}_{s_1 + 4 + j, i'}$ for all $i'$. Combining this with equations \eqref{ePreImp}, \eqref{eSelReminder} and \eqref{exReminder}, we have now shown \eqref{ePostImp}.
\end{proof}

\subsubsection{Auxiliary primes}
So far we have executed the first and second part of our construction of $\kappa$. We will now embark on the third part of our construction, where we add certain auxiliary primes. Pick a second subset $\{\tau_7, \ldots, \tau_{12}\}$ of $\{\sigma_1, \ldots, \sigma_{24}\}$ of cardinality $6$, disjoint from $\{\tau_1, \ldots, \tau_6\}$, such that
\begin{align}
\label{eRealDuplication}
\tau_i(x) \cdot \tau_{6 + i}(x) > 0
\end{align}
for all $x \in \{-1, \alpha, \beta, \gamma\}$ and for all $i \in \{1, \ldots, 6\}$. This is indeed possible; for $x = -1$ equation \eqref{eRealDuplication} always holds, while for $x \in \{\alpha, \beta, \gamma\}$, we use that $E$ satisfies Lemma \ref{lStartingCurve}$(b)$.

\begin{mydef}
Let $r \geq 0$ and let $\mathfrak{p}_1, \dots, \mathfrak{p}_r \not \in T$. We say that $p_{r + 1}, \dots, p_{r + 6}$ are auxiliary primes if
\begin{enumerate}
\item[$(C1)$] the ideals $\mathfrak{p}_{r + 1} := (p_{r + 1}), \dots, \mathfrak{p}_{r + 6} := (p_{r + 6})$ are distinct prime ideals coprime to $T$ and $\mathfrak{p}_1, \dots, \mathfrak{p}_r$;
\item[$(C2)$] we have for all $i \in \{1, \dots, 6\}$
\begin{alignat*}{3}
&\tau_i(p_{r + 1}) \cdot \tau_{6 + i}(p_{r + 1}) > 0, \quad \quad \quad \quad && \tau_i(p_{r + 1}) < 0 &&\Longleftrightarrow i = 1, \\
&\tau_i(p_{r + 2}) \cdot \tau_{6 + i}(p_{r + 2}) > 0, \quad \quad \quad \quad && \tau_i(p_{r + 2}) < 0 &&\Longleftrightarrow i = 2, \\
&\tau_i(p_{r + 3}) \cdot \tau_{6 + i}(p_{r + 3}) > 0, \quad \quad \quad \quad && \tau_i(p_{r + 3}) < 0 &&\Longleftrightarrow i = 3, \\
&\tau_i(p_{r + 4}) \cdot \tau_{6 + i}(p_{r + 4}) > 0, \quad \quad \quad \quad && \tau_i(p_{r + 4}) < 0 &&\Longleftrightarrow i = 4, \\
&\tau_i(p_{r + 5}) \cdot \tau_{6 + i}(p_{r + 5}) > 0, \quad \quad \quad \quad && \tau_i(p_{r + 5}) < 0 &&\Longleftrightarrow i = 5, \\
&\tau_i(p_{r + 6}) \cdot \tau_{6 + i}(p_{r + 6}) > 0, \quad \quad \quad \quad && \tau_i(p_{r + 6}) < 0 &&\Longleftrightarrow i = 6;
\end{alignat*}
\item[$(C3)$] for all real places $\tau \not \in \{\tau_1, \dots, \tau_{12}\}$ and for every $i \in \{1, \dots, 6\}$ we have $\tau(p_{r + i}) > 0$;
\item[$(C4)$] we have $p_{r + 1}, \ldots, p_{r + 6} \in K_v^{\ast 2}$ for all $v \in \{\mathfrak{p}_1, \ldots, \mathfrak{p}_r\}$ and for all finite places $v \in T$. Furthermore, for each $1 \leq i < j \leq 6$ we have
$$
\left(\frac{p_{r + j}}{p_{r + i}}\right) = -1 \Longleftrightarrow (i, j) \in \{(1, 2), (3, 4), (5, 6)\}.
$$
\end{enumerate}
\end{mydef}

\noindent Given $r \geq 0$ and distinct prime ideals $\mathfrak{p}_1, \dots, \mathfrak{p}_r \not \in T$, repeated application of Mitsui's prime ideal theorem \cite[Theorem 1.2.1]{KaiM} shows that there always exist auxiliary primes $p_{r + 1}, \dots, p_{r + 6}$. We also pick
\begin{enumerate}
\item[$(C5)$] local uniformizers $\pi_{r + 1}, \dots, \pi_{r + 6}$ such that for all $i \in \{1, \dots, 6\}$
$$
\pi_{r + i} = p_{r + i} \epsilon_{r + i},
$$
where $\epsilon_{r + i}$ is a non-square unit.
\end{enumerate}

\noindent We now collect some standard properties of the Hilbert symbol.

\begin{lemma}
\label{lHilbert}
Let $v$ be an odd place and let $\pi \in K_v^\ast$ be a local uniformizer. Then:
\begin{enumerate}
\item[(a)] $(\epsilon_1, \epsilon_2)_v = 1$ for all $\epsilon_1, \epsilon_2 \in K_v^\ast$ with even valuation;
\item[(b)] for all $\epsilon \in K_v^\ast$ with even valuation, we have $(\epsilon, \pi)_v = 1$ if and only if $\epsilon \in K_v^{\ast 2}$.
\end{enumerate}
Moreover, we have for all real places $w$ of $K$ and for all $x, y \in K_w^\ast$ 
$$
(x, y)_w = 1 \Longleftrightarrow x \textup{ and } y \textup{ are both negative.}
$$
\end{lemma}

\noindent We will start with some general facts about auxiliary primes. It is in this lemma that the places $\{\tau_7, \dots, \tau_{12}\}$ and equation \eqref{eRealDuplication} play a key role.

\begin{lemma}
\label{lAuxPrimes}
Let $p_{r + 1}, \dots, p_{r + 6}$ be auxiliary primes.
\begin{itemize}
\item[(i)] For all $i \in \{1, \dots, 6\}$, the elements $-1$, $\alpha$, $\beta$ and $\gamma$ are in $K_{p_{r + i}}^{\ast 2}$.
\item[(ii)] We have
$$
\left(\frac{p_{r + i}}{p_{r + j}}\right) = \left(\frac{p_{r + j}}{p_{r + i}}\right) \quad \quad \textup{ for all distinct } i, j \in \{1, \dots, 6\}.
$$
\item[(iii)] We have for all $i \in \{1, \dots, 6\}$
\begin{align}
\label{eFinalLocalSelmer2}
\mathcal{L}_{r + i, \boldsymbol{\pi}, \mathfrak{p}_{r + i}} = \langle (1, \epsilon_{r + i} p_{r + i}), (\epsilon_{r + i} p_{r + i}, 1) \rangle.
\end{align}
\end{itemize}
\end{lemma}

\begin{proof}
For $(i)$, we apply Hilbert reciprocity to the symbol $(p_{r + i}, x)$ for each $i \in \{1, \dots, 6\}$ and $x \in \{-1, \alpha, \beta, \gamma \}$. Note that $-1, \alpha, \beta, \gamma$ are units at all places outside of $T$. Hence by $(C4)$ and Lemma \ref{lHilbert}$(a)$, the symbol $(p_{r + i}, x)_v$ certainly vanishes at all finite places $v$ different from $p_{r + i}$. Moreover, the total contribution from the infinite places is trivial because of equation \eqref{eRealDuplication}, $(C2)$ and $(C3)$. In all cases, Hilbert reciprocity gives the vanishing of the local Hilbert symbol $(p_{r + i}, x)_{(p_{r + i})}$. Since $x$ has $p_{r + i}$-adic valuation zero, this gives part $(i)$ thanks to Lemma \ref{lHilbert}$(b)$.

Part $(ii)$ follows from Hilbert reciprocity applied to the symbol $(p_{r + i}, p_{r + j})$. Note that this symbol is trivial at all $v \in \Omega_K - \{p_{r + i}, p_{r + j}\}$; indeed, we have 
\begin{itemize}
\item $(p_{r + i}, p_{r + j})_v = 1$ for all odd places $v \in \Omega_K - \{p_{r + i}, p_{r + j}\}$ by Lemma \ref{lAuxPrimes}$(a)$;
\item $(p_{r + i}, p_{r + j})_v = 1$ for all $2$-adic places $v$ by the first part of $(C4)$;
\item $(p_{r + i}, p_{r + j})_v = 1$ for all infinite places $v$, as $p_{r + i}$ and $p_{r + j}$ are not simultaneously negative at $v$ by $(C2)$ and $(C3)$ (hence the result follows from the last part of Lemma \ref{lHilbert}).
\end{itemize}
Thus, part $(ii)$ follows from Hilbert reciprocity and Lemma \ref{lHilbert}$(b)$. Part $(iii)$ is a consequence of $(i)$, condition $(C5)$ and equation \eqref{eSelmerStructure}.
\end{proof}

Given auxiliary primes $p_{r + 1}, \dots, p_{r + 6}$, we now write
$$
V := \mathrm{Sel}_{\mathcal{L}_{r + 6, \boldsymbol{\pi}}}(G_K, E[2]).
$$
If we take two elements $w_1, w_2 \in V$, then these are always of the shape
\begin{align}
\label{eExpandab}
w_1 = \left(a \prod_{i = 1}^6 p_{r + i}^{e_{i, 1}}, b \prod_{i = 1}^6 p_{r + i}^{f_{i, 1}}\right), \quad \quad w_2 = \left(c \prod_{i = 1}^6 p_{r + i}^{e_{i, 2}}, d \prod_{i = 1}^6 p_{r + i}^{f_{i, 2}}\right)
\end{align}
with $(a, b)$ and $(c, d)$ having even valuation at all places $v \not \in T \cup \{\mathfrak{p}_1, \dots, \mathfrak{p}_r\}$.

\begin{lemma}
\label{c1}
Let $r \geq 0$, let $\mathfrak{p}_1, \dots, \mathfrak{p}_r \not \in T$ be distinct prime ideals and let $\pi_1, \dots, \pi_r$ be local uniformizers. Assume that
\begin{equation}
\label{eDiagonalGone2}
\mathcal{D}_{r ,i} = 1 \quad \textup{ for all } i \in \{1, 2, 3\}.
\end{equation}
Let $p_{r + 1}, \dots, p_{r + 6}$ be auxiliary primes. Let $w_1 \in V$ be non-trivial and write it as in equation \eqref{eExpandab}. Then $(a, b)$ is not in the kernel of the maps $T_1$, $T_2$ or $T_1 + T_2$ from Lemma \ref{lReduction}. Moreover, if $w_2 \in V$ does not equal $w_1$, then $(a, b)$ and $(c, d)$ are distinct (with $(c, d)$ as in equation \eqref{eExpandab}).
\end{lemma}

\begin{proof}
We start by showing that for all $e_1, \dots, e_6, f_1, \dots, f_6 \in \mathbb{F}_2$
\begin{equation}
\label{eVImpl}
\left(\prod_{i = 1}^6 p_{r + i}^{e_i}, \prod_{i = 1}^6 p_{r + i}^{f_i}\right) \in V \Longrightarrow e_i = f_i = 0.
\end{equation}
Using the local conditions at $\tau_1, \dots, \tau_6$ (see the top row of the table \eqref{eInfMat1}) and the sign assignments in $(C2)$, we deduce that $f_1 = f_3 = f_5 = 0$ and $e_2 = e_4 = e_6 = 0$. We will now show that $f_2 = 0$. Lemma \ref{lAuxPrimes}$(iii)$ and $f_1 = f_3 = f_5 = 0$ imply
\begin{equation}
\label{eLocals14}
\prod_{i = 1}^6 \res_{\mathfrak{p}_{r + 2}}(p_{r + i}^{f_i}) = \res_{\mathfrak{p}_{r + 2}}(p_{r + 2}^{f_2} p_{r + 4}^{f_4} p_{r + 6}^{f_6}) \in \langle \epsilon_{r + 2} \res_{\mathfrak{p}_{r + 2}}(p_{r + 2}) \rangle. 
\end{equation}
The second part of $(C4)$ reveals that $p_{r + 4}$ and $p_{r + 6}$ are in $K_{p_{r + 2}}^{\ast 2}$. Combining this with equation \eqref{eLocals14}, this proves $f_2 = 0$. Using Lemma \ref{lAuxPrimes}$(ii)$ if necessary, we may similarly prove that $f_4 = f_6 = 0$ and $e_1 = e_3 = e_5 = 0$, thus proving the implication \eqref{eVImpl}.

For the first part of Lemma \ref{c1}, let $w_1 \in V$ be non-trivial and write it as in equation \eqref{eExpandab}. Note that the auxiliary primes $p_{r + i}$ are squares at $\mathfrak{p}_1, \dots, \mathfrak{p}_r$ by the first part of $(C4)$. Hence $(a, b)$ lies in $\mathcal{L}_{r, \boldsymbol{\pi}, \mathfrak{p}_j}$ for all $j = 1, \dots, r$. Therefore, if $(a, b)$ were to be in the kernel of $T_1$, $T_2$ or $T_1 + T_2$, then $(a, b)$ must be trivial by equation \eqref{eDiagonalGone2}. Thus $w_1$ is trivial by \eqref{eVImpl}, contradicting our assumption.

For the second part of Lemma \ref{c1}, let $w_1, w_2 \in V$. We prove the contrapositive, i.e.~$(a, b) = (c, d)$ implies $w_1 = w_2$. So suppose that $(a, b) = (c, d)$. Looking at $w_1/w_2$ and applying \eqref{eVImpl}, we conclude that $w_1 = w_2$ as desired.
\end{proof}

We will now reap the benefits from Lemma \ref{c1}. In our next result we shall insist on adding exactly six primes; this is not strictly necessary but is convenient for numbering our primes in the last part of the argument.

\begin{lemma}
\label{lPrepDone}
Let $r \geq 0$, let $\mathfrak{p}_1, \dots, \mathfrak{p}_r \not \in T$ be distinct prime ideals and let $\pi_1, \dots, \pi_r$ be local uniformizers. Assume that
\begin{equation}
\label{eTrivialAtr}
\mathrm{Sel}_{\mathcal{L}_{r, \boldsymbol{\pi}}}(G_K, E[2]) = 1, \quad \quad \mathcal{D}_{r ,i} = 1 \quad \textup{ for all } i \in \{1, 2, 3\}.    
\end{equation}
Let $p_{r + 1}, \dots, p_{r + 6}$ be auxiliary primes. Then there exist distinct prime ideals $\mathfrak{p}_{r + 7}, \dots, \mathfrak{p}_{r + 12} \not \in T \cup \{\mathfrak{p}_1, \dots, \mathfrak{p}_r, (p_{r + 1}), \dots, (p_{r + 6})\}$ and local uniformizers $\pi_{r + 7}, \dots, \pi_{r + 12}$ such that
\begin{itemize}
\item we have
\begin{equation}
\label{eSelmerPrep}
\mathrm{Sel}_{\mathcal{L}_{r + 12, \boldsymbol{\pi}}}(G_K, E[2]) = 1,
\end{equation}
\item we have for all $i \in \{1, \dots, 6\}$ and all $j \in \{7, \dots, 12\}$
\begin{equation}
\label{eSelmerPrep2}
\left(\frac{p_{r + i}}{\mathfrak{p}_{r + j}}\right) = 1.
\end{equation}
\end{itemize}
\end{lemma}

\begin{proof}
We first show that for each $2$-dimensional subspace $W$ of $V = \mathrm{Sel}_{\mathcal{L}_{r + 6, \boldsymbol{\pi}}}(G_K, E[2])$, there exists a prime ideal $\mathfrak{P} \not \in T \cup \{\mathfrak{p}_1, \dots, \mathfrak{p}_r, (p_{r + 1}), \dots, (p_{r + 6})\}$ such that the auxiliary primes $p_{r + i}$ are squares modulo $\mathfrak{P}$ and such that the restriction map $\res_{\mathfrak{P}}|_W: W \rightarrow H^1(G_{K_{\mathfrak{P}}}, E[2])$ is injective.

To this end, take such a $2$-dimensional subspace $W$ of $V$ spanned by $w_1$ and $w_2$. Write $(a, b)$ and $(c, d)$ for the elements corresponding to $w_1$ and $w_2$ as in equation \eqref{eExpandab}. Let $K_1/K$ be the maximal abelian extension of $K$ of exponent $2$ that is unramified outside of $T \cup \{\mathfrak{p}_1, \dots, \mathfrak{p}_r\}$, and write
$$
K_2 := K(\{\sqrt{p_i} : i \in \{r + 1, \dots, r + 6\}\}).
$$
Then the extensions $K_1/K$ and $K_2/K$ are linearly disjoint by ramification considerations, and moreover 
$$
K(\sqrt{a}, \sqrt{b}, \sqrt{c}, \sqrt{d}) \subseteq K_1. 
$$
The assumption \eqref{eDiagonalGone2} of Lemma \ref{c1} is satisfied thanks to equation \eqref{eTrivialAtr}. By the first part of Lemma \ref{c1}, $\langle (a, b), (c, d) \rangle$ is not contained in the kernel of $T_1$, $T_2$ or $T_1 + T_2$. Moreover, by the second part of Lemma \ref{c1}, the dimension of $\langle (a, b), (c, d) \rangle$ equals $2$. Then, by the contrapositive of Lemma \ref{lReduction} applied to the quadratic characters $(\chi_1, \chi_2)$ and $(\chi_1', \chi_2')$ associated to respectively $(a, b)$ and $(c, d)$, we deduce the existence of a prime $\mathfrak{q} \not \in T \cup \{\mathfrak{p}_1, \dots, \mathfrak{p}_r, (p_{r + 1}), \dots, (p_{r + 6})\}$ such that
$$
\det
\begin{pmatrix}
\chi_1(\Frob_\mathfrak{q}) & \chi_2(\Frob_\mathfrak{q}) \\
\chi_1'(\Frob_\mathfrak{q}) & \chi_2'(\Frob_\mathfrak{q})
\end{pmatrix}
= 1.
$$
Recalling that the extensions $K_1/K$ and $K_2/K$ are linearly disjoint, the existence of $\mathfrak{P}$ follows from applying the Chebotarev density theorem to $K_1K_2/K$.

We will now prove the lemma. If $\dim_{\FF_2} V \geq 2$, then we fix a $2$-dimensional subspace $W \subseteq V$, and take $\mathfrak{p}_{r + 7}$ to be a prime $\mathfrak{P}$ as above. Since the restriction map $\res_{\mathfrak{P}}|_W: W \rightarrow H^1(G_{K_{\mathfrak{P}}}, E[2])$ is injective and since $W$ has dimension $2$, we have $n_{r + 6} = -2$ by Remark \ref{rSelmerChangeUpgrade}$(c)$ and Lemma \ref{lSelmerChange}. Moreover, given that $n_{r + 6} = -2$, it follows from Remark \ref{rSelmerChangeUpgrade}$(b)$ that
$$
\mathrm{Sel}_{\mathcal{L}_{r + 7, \boldsymbol{\pi}}}(G_K, E[2]) \subseteq \mathrm{Sel}_{\mathcal{L}_{r + 6, \boldsymbol{\pi}}}(G_K, E[2]) = V.
$$
Note that $\dim_{\FF_2} V$ is even by Lemma \ref{lParity}. Hence iterating this procedure, we find a sequence of prime ideals $\mathfrak{p}_{r + 7}, \dots, \mathfrak{p}_{r + 6 + \frac{1}{2}\dim_{\FF_2} V}$ such that
$$
\mathrm{Sel}_{\mathcal{L}_{r + 6 + \frac{1}{2}\dim_{\FF_2} V, \boldsymbol{\pi}}}(G_K, E[2]) = 1, \quad \quad \left(\frac{p_{r + i}}{\mathfrak{p}_{r + j}}\right) = 1 
$$
for all $i \in \{1, \dots, 6\}$ and all $j \in \{7, \dots, 6 + \frac{1}{2}\dim_{\FF_2} V\}$.

Note that we have the bound $\dim_{\FF_2} V \leq 12$ as a consequence of \eqref{eTrivialAtr} and Lemma \ref{lSelmerChange}. Set $k := s_2 + 6 + \frac{1}{2}\dim_{\FF_2} V \leq s_2 + 12$. If $k < s_2 + 12$, we pick $\mathfrak{p}_{k + 1}$ such that all $p_{r + i}$ (with $i \in \{1, \dots, 6\}$) are squares modulo $\mathfrak{p}_{k + 1}$. After choosing the local uniformizer $\pi_{k + 1}$ so that $\mathcal{L}_{k + 1, \boldsymbol{\pi}, \mathfrak{p}_{k + 1}} \neq A_{\mathfrak{p}_{k + 1}}$, we are not in the first case of Lemma \ref{lSelmerChange}. Hence if $\mathrm{Sel}_{\mathcal{L}_{k, \boldsymbol{\pi}}}(G_K, E[2]) = 1$, we must be in the third case and thus also $\mathrm{Sel}_{\mathcal{L}_{k + 1, \boldsymbol{\pi}}}(G_K, E[2]) = 1$. Iterating this procedure exactly $6 - \frac{1}{2} \dim_{\FF_2} V$ times gives the lemma.
\end{proof}

\subsubsection{The final step}
We are now ready to finish the construction of $\kappa$.

\begin{theorem}
There exists an auxiliary twist $\kappa$.
\end{theorem}

\begin{proof}
Fix an integer $s_2 \geq 0$, distinct prime ideals $\mathfrak{p}_1, \dots, \mathfrak{p}_{s_2} \not \in T$ and local uniformizers $\pi_1, \dots, \pi_{s_2}$ as in Lemma \ref{lAux2}. We will take $s := s_2 + 16$. Let $p_{s_2 + 1}, \dots, p_{s_2 + 6}$ be auxiliary primes and let $\pi_{s_2 + 1}, \dots, \pi_{s_2 + 6}$ be local uniformizers as in $(C5)$. Let $\mathfrak{p}_{s_2 + 7}, \dots, \mathfrak{p}_{s_2 + 12}$ and $\pi_{s_2 + 7}, \dots, \pi_{s_2 + 12}$ be respectively the primes and local uniformizers guaranteed by Lemma \ref{lPrepDone}. We now list a series of conditions for our next three prime elements $p_{s_2 + 13}, p_{s_2 + 14}, p_{s_2 + 15}$
\begin{enumerate}
\item[$(C6)$] the ideals $\mathfrak{p}_{s_2 + 13} := (p_{s_2 + 13})$, $\mathfrak{p}_{s_2 + 14} := (p_{s_2 + 14})$ and $\mathfrak{p}_{s_2 + 15} := (p_{s_2 + 15})$ are distinct prime ideals coprime to $T$ and $\mathfrak{p}_1, \dots, \mathfrak{p}_{s_2 + 12}$;
\item[$(C7)$] we have $p_{s_2 + 13}, p_{s_2 + 14}, p_{s_2 + 15} \in K_v^{\ast 2}$ for all $v \in T \cup \{\mathfrak{p}_1, \dots, \mathfrak{p}_{s_2}\} \cup \{\mathfrak{p}_{s_2 + 7}, \dots, \mathfrak{p}_{s_2 + 12}\}$;
\item[$(C8)$] for each $1 \leq i \leq 6$ and each $1 \leq j \leq 3$ we have that
$$
\left(\frac{p_{s_2 + 12 + j}}{p_{s_2 + i}}\right) = -1 \Longleftrightarrow (i, j) \in \left\{(1, 1), (2, 1), (3, 2), (4, 2), (5, 3), (6, 3)\right\};
$$
\item[$(C9)$] for each $1 \leq i < j \leq 3$, we have that
$$
\left(\frac{p_{s_2 + 12 + j}}{p_{s_2 + 12 + i}}\right) = 1.
$$
\end{enumerate}

\noindent It is clearly possible to find such prime elements by repeatedly applying Mitsui's prime ideal theorem \cite[Theorem 1.2.1]{KaiM}. We next pick
\begin{enumerate}
\item[$(C10)$] local uniformizers $\pi_{s_2 + 13}, \pi_{s_2 + 14}, \pi_{s_2 + 15}$ such that for all $i \in \{1, 2, 3\}$
$$
\pi_{s_2 + 12 + i} = p_{s_2 + 12 + i} \epsilon_{s_2 + 12 + i},
$$
where $\epsilon_{s_2 + 12 + i}$ is a non-square unit.
\end{enumerate}

\noindent We start by collecting some basic facts about the prime elements $p_{s_2 + 13}, \dots, p_{s_2 + 15}$. 

\begin{lemma}
\label{c2}
\begin{itemize}
\item[(i)] For all $1 \leq i \leq 6$ and all $1 \leq j \leq 3$
$$
\left(\frac{p_{s_2 + i}}{p_{s_2 + 12 + j}}\right) = \left(\frac{p_{s_2 + 12 + j}}{p_{s_2 + i}}\right).
$$
\item[(ii)] The elements $-1$, $\alpha$, $\beta$ and $\gamma$ lie in $K_{p_{s_2 + 13}}^{\ast 2}$, $K_{p_{s_2 + 14}}^{\ast 2}$ and $K_{p_{s_2 + 15}}^{\ast 2}$.
\item[(iii)] For all $i \in \{s_2 + 13, s_2 + 14, s_2 + 15\}$
$$
\mathcal{L}_{s_2 + 15, \boldsymbol{\pi}, \mathfrak{p}_i} = \langle (1, \epsilon_i p_i), (\epsilon_i p_i, 1) \rangle.
$$
\end{itemize} 
\end{lemma}

\begin{proof}
Our arguments will be rather similar to Lemma \ref{lAuxPrimes}. For $(i)$, we apply Hilbert reciprocity to the symbol $(p_{s_2 + i}, p_{s_2 + 12 + j})$. The contribution from odd places $v \in \Omega_K - \{(p_{s_2 + i}), (p_{s_2 + 12 + j})\}$ vanishes by Lemma \ref{lHilbert}$(a)$. Next, at the $2$-adic and infinite places, we also get vanishing of the local Hilbert symbol by $(C7)$. Hence $(i)$ follows from an examination of the local Hilbert symbol at $v \in \{(p_{s_2 + i}), (p_{s_2 + 12 + j})\}$ and an application of Lemma \ref{lHilbert}$(b)$.

For $(ii)$, we apply Hilbert reciprocity to the symbol $(p_{s_2 + 12 + j}, x)$ for each $j \in \{1, 2, 3\}$ and $x \in \{-1, \alpha, \beta, \gamma \}$. This symbol vanishes by $(C7)$ and Lemma \ref{lHilbert}$(a)$ at all places $v \in \Omega_K - \{(p_{s_2 + 12 + j})\}$. Then the result follows from Lemma \ref{lHilbert}$(b)$ applied to the place $v := (p_{s_2 + 12 + j})$. For $(iii)$, it follows from $(C10)$, $(ii)$ and the defining equation \eqref{eSelmerStructure} of $\mathcal{L}_{s_2 + 15, \boldsymbol{\pi}, \mathfrak{p}_i}$ that
$$
\mathcal{L}_{s_2 + 15, \boldsymbol{\pi}, \mathfrak{p}_i} = \langle (1, \epsilon_i p_i), (\epsilon_i p_i, 1) \rangle
$$
for all $i \in \{s_2 + 13, s_2 + 14, s_2 + 15\}$. 
\end{proof}

\noindent We claim that the six elements
\begin{alignat}{6}
&(z_1, z_2) &&:= (p_{s_2 + 1} p_{s_2 + 13}, 1), \ &&(z_3, z_4) &&:= (1, p_{s_2 + 2} p_{s_2 + 13}), \ &&(z_5, z_6) &&:= (p_{s_2 + 3} p_{s_2 + 14}, 1) \nonumber \\
&(z_7, z_8) &&:= (1, p_{s_2 + 4} p_{s_2 + 14}), \ &&(z_9, z_{10}) &&:= (p_{s_2 + 5} p_{s_2 + 15}, 1), \ &&(z_{11}, z_{12}) &&:= (1, p_{s_2 + 6} p_{s_2 + 15}) \label{eFinalBasis}
\end{alignat}
form a basis of $\mathrm{Sel}_{\mathcal{L}_{s_2 + 15, \boldsymbol{\pi}}}(G_K, E[2])$ and satisfy the sign conditions in the table \eqref{eInfMat1}. Note that these six elements are visibly linearly independent, and that we have the inequality
$$
\dim_{\mathbb{F}_2} \mathrm{Sel}_{\mathcal{L}_{s_2 + 15, \boldsymbol{\pi}}}(G_K, E[2]) \leq 6
$$
by Lemma \ref{lSelmerChange} and by equation \eqref{eSelmerPrep} (applied with $r = s_2$). So to establish the claim, it suffices to check that the above six elements restrict to $\mathcal{L}_{s_2 + 15, \boldsymbol{\pi}, v}$ for all $v \in \Omega_K$ and satisfy the sign conditions in the table \eqref{eInfMat1}.

For the sign conditions, these are a direct consequence of conditions $(C2)$ and $(C7)$. In particular, the above six elements restrict to $\mathcal{L}_{s_2 + 15, \boldsymbol{\pi}, \tau}$ for $\tau \in \{\tau_1, \dots, \tau_6\}$. Moreover, $(z_1, z_2), \dots, (z_{11}, z_{12})$ are then also in $\mathcal{L}_{s_2 + 15, \boldsymbol{\pi}, \tau}$ for $\tau \in \{\tau_7, \dots, \tau_{12}\}$ by equation \eqref{eRealDuplication} and conditions $(C2)$ and $(C7)$. Finally, we see that $(z_1, z_2), \dots, (z_{11}, z_{12})$ are in $\mathcal{L}_{s_2 + 15, \boldsymbol{\pi}, \tau}$ at all other archimedean places $\tau$ by $(C3)$ and $(C7)$.

Observe that the elements $z_1, \dots, z_{12}$ have even valuation at all places $v$ outside $T \cup \{\mathfrak{p}_1, \ldots, \mathfrak{p}_{s_2 + 15}\}$, and therefore the pairs $(z_1, z_2), \dots, (z_{11}, z_{12})$ are in $\mathcal{L}_{s_2 + 15, \boldsymbol{\pi}, v}$ for all such $v$ by definition. We now turn to the finite places in $T \cup \{\mathfrak{p}_1, \ldots, \mathfrak{p}_{s_2}\}$, which are covered by $(C4)$ and $(C7)$, so it remains to check that $(z_1, z_2), \dots, (z_{11}, z_{12})$ are in $\mathcal{L}_{s_2 + 15, \boldsymbol{\pi}, \mathfrak{p}}$ for all $\mathfrak{p} \in \{\mathfrak{p}_{s_2 + 1}, \dots, \mathfrak{p}_{s_2 + 15}\}$. For the prime ideals $\mathfrak{p}_{s_2 + 7}, \dots, \mathfrak{p}_{s_2 + 12}$, this is a consequence of equation \eqref{eSelmerPrep2} and condition $(C7)$. 

We now treat the remaining places $\mathfrak{p}_{s_2 + 1}, \dots, \mathfrak{p}_{s_2 + 6}, \mathfrak{p}_{s_2 + 13}, \mathfrak{p}_{s_2 + 14}, \mathfrak{p}_{s_2 + 15}$. By combining Lemma \ref{lAuxPrimes}$(iii)$ and Lemma \ref{c2}$(iii)$, we conclude that
$$
\mathcal{L}_{s_2 + 15, \boldsymbol{\pi}, \mathfrak{p}_i} = \langle 1, \epsilon_i p_i), (\epsilon_i p_i, 1) \rangle
$$
for all $i \in \{s_2 + 1, \dots, s_2 + 6, s_2 + 13, s_2 + 14, s_2 + 15\}$. With this proven, one checks directly that $\res_{\mathfrak{p}_i}(z_{2j + 1}, z_{2j + 2}) \in \mathcal{L}_{i, \boldsymbol{\pi}, \mathfrak{p}_i}$ by using Lemma \ref{c2}$(i)$, $(C4)$, $(C8)$ and $(C9)$. This ends the proof that the elements in equation \eqref{eFinalBasis} form a basis of $\mathrm{Sel}_{\mathcal{L}_{s_2 + 15, \boldsymbol{\pi}}}(G_K, E[2])$ and satisfy the sign conditions in the table \eqref{eInfMat1}. 

Let $\mathfrak{m}$ be the modulus $8 \infty \mathfrak{r}$, where $\mathfrak{r}$ denotes the product of odd prime ideals in $T$. Denote by $K(\mathfrak{m})/K$ the corresponding ray class field and denote by
$$
L := K(\{\sqrt{p_i} : i \in \{s_2 + 1, \ldots, s_2 + 6\} \cup \{s_2 + 13, s_2 + 14, s_2 + 15\}\}).
$$
Then $K(\mathfrak{m})/K$ and $L/K$ are linearly disjoint by ramification considerations. Therefore, by the Chebotarev density theorem applied to $K(\mathfrak{m})L/K$, there exists a prime ideal $\mathfrak{p}_{s_2 + 16} \not \in T \cup \{\mathfrak{p}_1, \dots, \mathfrak{p}_{s_2 + 15}\}$ such that the ideal $\mathfrak{p}_1 \cdots \mathfrak{p}_{s_2 + 16}$ admits a generator $\kappa$ that is $1$ modulo $8$, is congruent to $1$ modulo every odd prime in $T$, is negative at $\tau_4, \tau_5, \tau_6$, is positive at all other real places, and $\mathfrak{p}_{s_2 + 16}$ splits completely in $L/K$. This latter condition is equivalent to all $p_i$ with $i \in \{s_2 + 1, \ldots, s_2 + 6\} \cup \{s_2 + 13, s_2 + 14, s_2 + 15\}$ being squares modulo $\mathfrak{p}_{s_2 + 16}$. Inspecting the elements in \eqref{eFinalBasis}, we see that they all reduce trivially modulo $\mathfrak{p}_{s_2 + 16}$. This implies the inclusion
\begin{equation}
\label{eInclSel}
\mathrm{Sel}_{\mathcal{L}_{s_2 + 15, \boldsymbol{\pi}}}(G_K, E[2]) \subseteq \mathrm{Sel}_{\mathcal{L}_{s_2 + 16, \boldsymbol{\pi}}}(G_K, E[2]).
\end{equation}
Hence we are either in the first or third case of Lemma \ref{lSelmerChange}. Choosing the local uniformizer $\pi_{s_2 + 16}$ such that $\mathcal{L}_{s_2 + 16, \boldsymbol{\pi}, \mathfrak{p}_{s_2 + 16}} \neq A_{\mathfrak{p}_{s_2 + 16}}$, we must be in the third case of Lemma \ref{lSelmerChange} and therefore $n_{s_2 + 16} = 0$. Thus we deduce from equation \eqref{eInclSel} that
$$
\mathrm{Sel}_{\mathcal{L}_{s_2 + 16, \boldsymbol{\pi}}}(G_K, E[2]) = \mathrm{Sel}_{\mathcal{L}_{s_2 + 15, \boldsymbol{\pi}}}(G_K, E[2]).
$$
In particular, the elements in \eqref{eFinalBasis} form a basis of $\mathrm{Sel}_{\mathcal{L}_{s_2 + 16, \boldsymbol{\pi}}}(G_K, E[2])$ as well. Then $\kappa$ satisfies $(K1)$ and $(K2)$, ending the proof.
\end{proof}

\appendix
\newpage

\section{Additive combinatorics}
In this appendix, we shall first introduce the necessary notation to formally state Kai's result \cite[Theorem 12.1]{KaiAppendix} in the special case of four linear forms. The additive combinatorics machinery of Kai builds on the earlier works \cite{GT0A, GT1A, GT2A, GT3A}. After stating Kai's result, we shall then derive the additive combinatorics result that we shall use in this paper. We start by introducing two classes of norms that will be important to us.

\begin{mydef}
Let $K$ be a number field. We define the canonical norm on $K$ to be
$$
||x||_{\textup{can}} := \max_{\sigma: K \xhookrightarrow{} \C} |\sigma(x)|
$$
for $x \in K$. Here the maximum ranges over all embeddings $K \xhookrightarrow{} \C$.
\end{mydef}

\begin{mydef}
Let $K$ be a number field and let $\mathfrak{a}$ be a non-zero fractional ideal of $K$. Fix an ordered $\Z$-basis $\omega := (\omega_1, \dots, \omega_n)$ of $\mathfrak{a}$. This gives rise to a norm on $\mathfrak{a}$ by defining the norm of $\alpha = x_1 \omega_1 + \dots + x_n \omega_n \in \mathfrak{a}$ (with $x_1, \dots, x_n \in \Z$) to be
$$
||\alpha||_{\omega ,\infty} = \max(|x_1|, \dots, |x_n|).
$$
\end{mydef}

Since Kai is interested in uniformity over the fractional ideals $\mathfrak{a}$, the notion of norm-length compatible basis is introduced. The point of this notion is that one wants to uniformly pick bases of $\mathfrak{a}$ that are close to the canonical norm. For us, this uniformity is not relevant as we will ultimately take $\mathfrak{a} = O_K$. Nevertheless, we shall need to introduce this notion to formally state his result.

\begin{mydef}
Let $C > 1$. Then an ordered $\Z$-basis $\omega := (\omega_1, \dots, \omega_n)$ of $\mathfrak{a}$ is said to be $C$-norm-length compatible if
$$
C^{-1} N_{K/\Q}(\mathfrak{a})^{1/n} ||x||_{\omega, \infty} \leq ||x||_{\textup{can}} \leq C N_{K/\Q}(\mathfrak{a})^{1/n} ||x||_{\omega, \infty}
$$
for all $x \in \mathfrak{a}$.
\end{mydef}

Kai shows that norm-length compatible bases always exist.

\begin{proposition}[{\cite[Proposition 3.1]{KaiAppendix}}]
\label{pComp}
Let $K$ be a number field. Then there is a constant $C_K > 1$ such that for all non-zero fractional ideals $\mathfrak{a}$ of $K$, there is a $C_K$-norm-length compatible basis of $\mathfrak{a}$.
\end{proposition}

Kai then proceeds to fix such a choice of $C_K$ throughout his paper, and always takes norm-length compatible bases of $\mathfrak{a}$. The advantage of doing this is that Kai's result holds uniformly over the choice of $\mathfrak{a}$, a flexibility that will not be of relevance for us. We now introduce two direct number field analogues of respectively the von Mangoldt and Euler phi function.

\begin{mydef}
For a number field $K$ and a non-zero integral ideal $I$ of $O_K$, we define
$$
\varphi(I) := |(O_K/I)^\ast|.
$$
We also define the von Mangoldt function on all integral ideals $I$ of $O_K$ through
$$
\Lambda(I) = 
\begin{cases}
\log N_{K/\Q}(\mathfrak{p}) &\textup{if } I = \mathfrak{p}^m \textup{ for some prime } \mathfrak{p} \textup{ and some } m \geq 1 \\
0 &\textup{otherwise.}
\end{cases}
$$
Given a fractional ideal $\mathfrak{a}$, we define
$$
\Lambda^\mathfrak{a}(x) = \Lambda(x \mathfrak{a}^{-1})
$$
for $x \in \mathfrak{a}$.
\end{mydef}

An affine linear form $\psi: \Z^d \rightarrow \Z^n$ can always be written as
$$
\psi(x_1, \dots, x_d) = c_0 + \sum_{i = 1}^d c_i x_i
$$
for some $c_0, \dots, c_d \in \Z^n$. We then define
$$
\psi^{\text{hom}}(x_1, \dots, x_d) := \sum_{i = 1}^d c_i x_i
$$
to be the corresponding homogeneous linear form and we define for each integer $N \geq 1$ another norm
$$
||\psi||_N = \max(||c_0||_\infty/N, ||c_1||_\infty, \dots, ||c_d||_\infty),
$$
where $||\cdot||_\infty$ denotes the usual infinity norm. We define
$$
\beta_p^{\mathfrak{a}}(\psi_1, \dots, \psi_4) := \frac{p^{4n}}{\varphi(p)^4} \frac{|\{x \in (\Z/p\Z)^d : \psi_i(x) \neq 0 \textup{ in } \mathfrak{a}/\mathfrak{p} \mathfrak{a} \textup{ for all } i \in \{1, \dots, 4\} \textup{ and } \mathfrak{p} \mid p\}|}{p^d}.
$$
Kai \cite[Remark 12.5]{KaiAppendix} proves that $\prod_p \beta_p^{\mathfrak{a}}(\psi_1, \dots, \psi_4)$ converges absolutely.

\begin{theorem}[{\cite[Theorem 12.1]{KaiAppendix}}]
\label{tKai}
For all $n \in \Z_{\geq 1}$, there is a positive real number $c_n > 0$ such that for all $d \in \Z_{\geq 1}$, there is a positive real number $c_{n, d} > 0$ such that for all number fields $K$ of degree $n$, there exists $C > 0$ such that the following holds. 

Let $\mathfrak{a}$ be a non-zero fractional ideal equipped with a $C_K$-norm-length compatible basis. Let $N > C$. Let $\psi_1, \dots, \psi_4: \Z^d \rightarrow \mathfrak{a} \cong \Z^n$ be affine linear forms, where the isomorphism $\mathfrak{a} \cong \Z^n$ comes from the norm-length compatible basis. We assume that:
\begin{itemize}
\item the restriction of $\psi_i^{\textup{hom}}$ to the kernel of $\psi_j^{\textup{hom}}$ has finite cokernel for all $i \neq j$, and
\item $||\psi_i||_N \leq (\log \log N)^{\frac{1}{10^3 d^2} c_n}$.
\end{itemize} 
Let $\Omega \subseteq [-N, N]^d$ be a convex set. Then we have
$$
\left|\sum_{x \in \Omega \cap \Z^d} \left(\prod_{i = 1}^4 \Lambda^{\mathfrak{a}}(\psi_i(x))\right) - \frac{\Vol(\Omega)}{\textup{res}_{s = 1} \zeta_K(s)^4} \cdot \prod_p \beta_p^{\mathfrak{a}}(\psi_1, \dots, \psi_4)\right| \leq \frac{CN^d}{(\log \log N)^{c_{n, d}}}.
$$
\end{theorem}

Since our affine linear forms will be fixed and since we will take $\mathfrak{a} = O_K$, let us state a version of this result without any uniformity.

\begin{theorem}
\label{tKai2}
Let $K$ be a number field and let $d \geq 1$. Let $\psi_1, \dots, \psi_4: \Z^d \rightarrow O_K$ be affine linear forms. We assume that the restriction of $\psi_i^{\textup{hom}}$ to the kernel of $\psi_j^{\textup{hom}}$ has finite cokernel for all $i \neq j$. Let $\Omega \subseteq \mathbb{R}^d$ be a convex set, whose volume grows as 
$$
\Vol(\Omega \cap [-N, N]^d) = cN^d + o(N^d)
$$
for some fixed $c > 0$ as $N \rightarrow \infty$. Then we have
$$
\sum_{x \in \Omega \cap [-N, N]^d \cap \Z^d} \left(\prod_{i = 1}^4 \Lambda^{O_K}(\psi_i(x))\right) = \frac{cN^d}{\textup{res}_{s = 1} \zeta_K(s)^4} \cdot \prod_p \beta_p^{O_K}(\psi_1, \dots, \psi_4) + o(N^d).
$$
\end{theorem}

\begin{proof}
This follows from Theorem \ref{tKai} and the fact that $O_K$ has a $C_K$-norm-length compatible basis by Proposition \ref{pComp}.
\end{proof}

\end{document}